\documentclass[journal]{IEEEtran}

\usepackage{graphicx}
\usepackage{pdfpages}
\ifCLASSINFOpdf
\else
\fi
\usepackage{float}
\usepackage{subfig}
\usepackage{epsfig}
\usepackage{amsmath}
\usepackage{mathtools}
\usepackage{amsfonts}
\usepackage{amssymb}
\usepackage{setspace}
\usepackage{tabularx}
\usepackage{algorithm,algorithmic}
\usepackage{enumerate}
\usepackage{dblfloatfix}
\usepackage{verbatim}
\begin{document}

\title{Accuracy Controlled Structure-Preserving ${\cal H}^2$-Matrix-Matrix Product in Linear Complexity with Change of Cluster Bases}

\author{Miaomiao~Ma,~\IEEEmembership{Graduate~Student~Member,~IEEE,}
        and~Dan~Jiao,~\IEEEmembership{Fellow,~IEEE}% <-this % stops a space
\thanks{The authors are with the School of
 Electrical and Computer Engineering, Purdue University, West Lafayette, IN, 47907 USA (e-mail: djiao@purdue.edu).}% <-this % stops a space
\thanks{Manuscript received July 31, 2019. This paper is an expanded version from the IEEE MTT-S International Conference on Numerical Electromagnetic and Multiphysics Modeling and Optimization, Cambridge, MA, USA, May 29-31, 2019.}
}%; revised August 26, 2018.

\maketitle

\begin{abstract}
${\cal H}^2$-matrix constitutes a general mathematical framework for efficient computation of both partial-differential-equation and integral-equation-based operators. Existing linear-complexity ${\cal H}^2$ matrix-matrix product (MMP) algorithm lacks explicit accuracy control, while controlling accuracy without compromising linear complexity is challenging. In this paper, we develop an accuracy controlled ${\cal H}^2$ matrix-matrix product algorithm by instantaneously changing the cluster bases during the matrix product computation based on prescribed accuracy. Meanwhile, we retain the computational complexity of the overall algorithm to be linear. Different from the existing ${\cal H}^2$ matrix-matrix product algorithm where formatted multiplications are performed using the original cluster bases, in the proposed algorithm, all additions and multiplications are either exact or computed based on prescribed accuracy. Furthermore, the original ${\cal H}^2$-matrix structure is preserved in the matrix product. While achieving optimal complexity for constant-rank matrices, the computational complexity of the proposed algorithm is also minimized for variable-rank ${\cal H}^2$-matrices. The proposed work serves as a fundamental arithmetic in the development of fast solvers for large-scale electromagnetic analysis. Applications to both large-scale capacitance extraction and electromagnetic scattering problems involving millions of unknowns on a single core have demonstrated the accuracy and efficiency of the proposed algorithm.
\end{abstract}


% \fi
\IEEEpeerreviewmaketitle

\section{Introduction}
\IEEEPARstart{T}{he} ${\cal H}^2$-matrix \cite{BormH2MMP,BormH2book} constitutes a general mathematical framework for compact representation and efficient computation of large dense systems. Both partial differential equation (PDE) and integral equation (IE) operators in electromagnetics can be represented as ${\cal H}^2$-matrices with controlled accuracy \cite{HaixinH, ZhouH, JiaoRank}.

The development of ${\cal H}^2$-matrix arithmetic such as addition, multiplication, and inverse are of critical importance to the development of fast solvers in electromagnetics \cite{Weninverse}. Under the ${\cal H}^2$-matrix framework, it has been shown that an ${\cal H}^2$-matrix-based addition, matrix-vector product (MVP), and matrix-matrix product (MMP) all can be performed in linear complexity for constant-rank ${\cal H}^2$ \cite{BormH2MMP}. However, the accuracy of existing ${\cal H}^2$-MMP algorithm like \cite{BormH2MMP} is not controlled. This is because given two ${\cal H}^2$-matrices $\textbf{A}_{{\cal H}^2}$ and $\textbf{B}_{{\cal H}^2}$, the matrix structure and cluster bases of their product $\textbf{C}=\textbf{A}_{{\cal H}^2} \times \textbf{B}_{{\cal H}^2}$ are pre-assumed, and a formatted multiplication is performed, whose accuracy is not controlled. For example, the row cluster bases of $\textbf{A}_{{\cal H}^2}$ and the column cluster bases of $\textbf{B}_{{\cal H}^2}$ are assumed to be those of $\textbf{C}$. This treatment lacks accuracy control since the original cluster basis may not be able to represent the new matrix content generated during the MMP. For instance, when multiplying a full-matrix block $\textbf{F}$ by a low rank block $\textbf{V}_t\textbf{S}{\textbf{V}_s}^T$, treating the result as a low-rank block is correct. However, it is inaccurate to use the original row cluster basis $\textbf{V}_t$ as the product's row cluster basis, since the latter has been changed to $\textbf{F}\textbf{V}_t$. Therefore, the algorithm in \cite{BormH2MMP} can be accurate if the cluster bases of the original matrices can also be used to accurately represent the matrix product. However, this is unknown in general applications, and hence the accuracy of existing linear-complexity MMP algorithm is not controlled. One can find many cases where a formatted multiplication would fail.

The posteriori multiplication in \cite{BormH2book} is more accurate than the formatted multiplication in \cite{BormH2MMP}. But it is only suitable for special ${\cal H}^2$ matrices. Besides, this posteriori multiplication requires much more computational time and memory than the formatted one. It needs to first present the product in ${\cal H}$-matrix and then convert it into an ${\cal H}^2$-matrix, the complexity of which is not linear.  %We can never encounter the situation that the target block $\textbf{C}_{i,k}$ is admissible, but the $\textbf{A}_{i,j}$ and $\textbf{B}_{j,k}$ are inadmissible. 

%Recently, in \cite{MiaoAPSH2MMP, MiaoICEAAH2MMP}, a new linear complexity algorithm is developed to compute the product of two  ${\cal H}^2$-matrices in controlled accuracy. The original cluster bases are changed by appending new cluster bases to account for the updates to the original matrix generated during the MMP. The product structure can be preserved as in \cite{MiaoICEAAH2MMP}, or can be changed based on the structure of the two multipliers as in \cite{MiaoAPSH2MMP}. However, the rank of the changed cluster bases by appending new ones can only be higher than the original rank.

In this work, we propose a new algorithm to do the ${\cal H}^2$ matrix-matrix multiplication with controlled accuracy. The cluster bases are calculated instantaneously based on the prescribed accuracy during the computation of the matrix-matrix product. Meanwhile, we are able to keep the computational complexity to be linear for constant-rank ${\cal H}^2$. For variable-rank cases such as those in an electrically large analysis, the proposed MMP is also efficient since it only involves $O(2^l)$ computations at level $l$, each of which costs $O(k_l^3)$ only, where $k_l$ is the rank at tree level $l$. This algorithm can be used as a fundamental arithmetic in the error-controlled fast inverse, LU factorization, solution for many right hand sides, etc. Numerical experiments have demonstrated its accuracy and low complexity. In \cite{MiaoAPSH2MMP, MiaoICEAAH2MMP}, we present a fast algorithm to compute the product of two ${\cal H}^2$-matrices in controlled accuracy. However, unlike this work, the original cluster bases are not completely changed, but appended to account for the updates to the original matrix during the MMP. In \cite{NEMO19}, we present the basic idea of this work. However, it is a one-page abstract. In this paper, we provide a complete algorithm together with a comprehensive analysis of its accuracy and complexity, whose validity and performance are then demonstrated by abundant numerical examples.

\section{PRELIMINARIES}
\begin{figure}[t]
  \centering
  \subfloat[]{\includegraphics[width=0.45\textwidth]{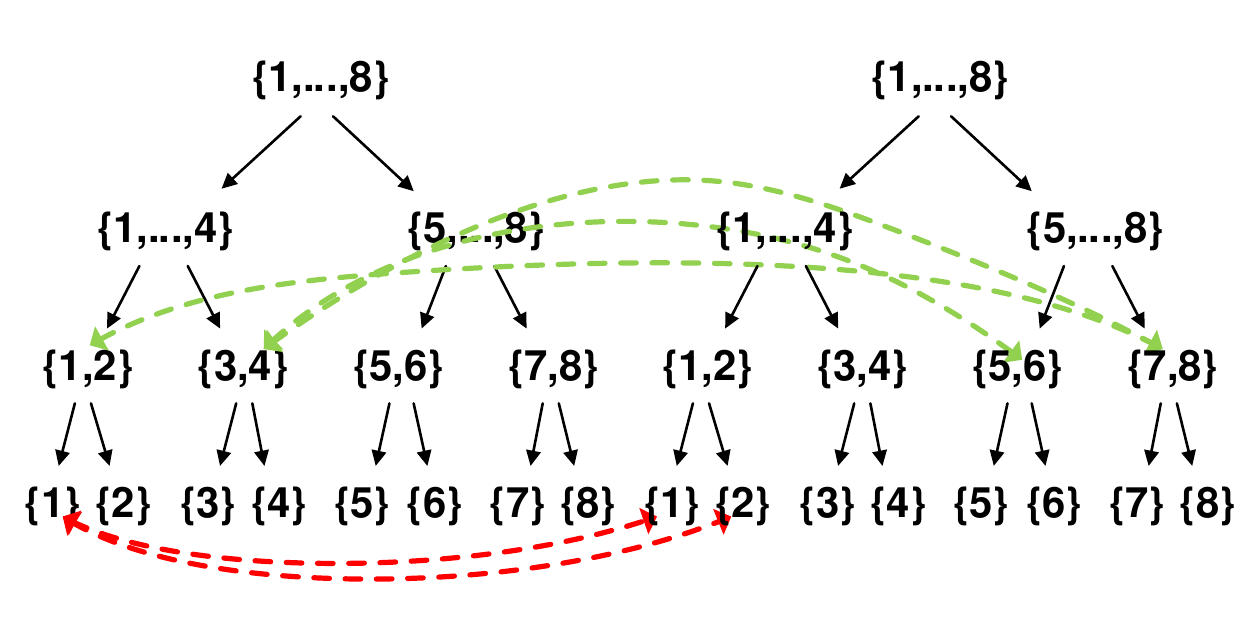}}
  \hskip 0.5truein
  \subfloat[]{\includegraphics[width=0.45\textwidth]{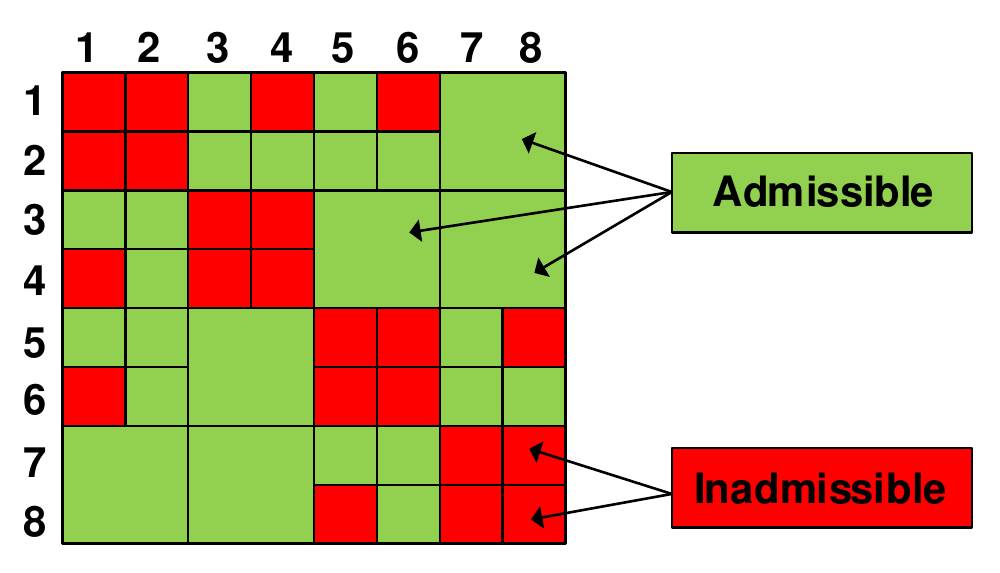}}
  \caption{Illustration of a block cluster tree and resulting ${\cal H}^2$-matrix partition. (a) Block cluster tree. (b) ${\cal H}^2$-matrix structure.}
\end{figure}
In an ${\cal H}^2$-matrix \cite{BormH2MMP}, the entire matrix is partitioned into multilevel admissible and inadmissible blocks, where inadmissible blocks are at the leaf level, noted as $\textbf{F}_{t,s}$. An admissible matrix block $\textbf{R}_{t,s}$ satisfies the following strong admissibility condition
\begin{equation}\label{admCond}
max\{diam( \Omega_t ), diam( \Omega_s ) \} \leq \eta dist(\Omega_t,\Omega_s),
\end{equation} 
where $\Omega_t$ ($\Omega_s$) denotes the geometrical support of the unknown set $t$ ($s$), $diam\{ \cdot \}$ is the Euclidean diameter of a set, $dist\{\cdot ,\cdot\}$ denotes the Euclidean distance between two sets, and $\eta$ is a positive parameter that can be used to control the admissibility condition. An admissible matrix block in an ${\cal H}^2$-matrix is represented as 
\begin{equation}
\textbf{R}_{t,s}=(\textbf{V}_t)_{\#t\times k} (\textbf{S}_{t,s})_{k\times k}(\textbf{V}_s)^{T}_{\#s\times k}\\
\end{equation}
where $\textbf{V}_{t}$ ($\textbf{V}_{s}$) is called cluster basis associated with cluster $t$ ($s$), $\textbf{S}_{t,s}$ is called coupling matrix. The cluster bases $\textbf{V}$ in an ${\cal H}^2$-matrix has a nested property. This means the cluster basis for a non-leaf cluster $t$, $\textbf{V}_{t}$, can be expressed by its two children's cluster bases, $\textbf{V}_{t_1}$ and $\textbf{V}_{t_2}$, as
\begin{equation}
(\textbf{V}_{t})_{\#t\times k} =
 \begin{bmatrix}
(\textbf{V}_{t_1})_{\#t_{1}\times k_{1}} & 0 \\
0 & (\textbf{V}_{t_2})_{\#t_{2}\times k_{2}} \\
\end{bmatrix} 
\begin{bmatrix}
(\textbf{T}_{t_1})_{k_{1} \times k} \\
(\textbf{T}_{t_2})_{k_{2} \times k} \\
\end{bmatrix}
\end{equation}
where $\textbf{T}_{t_1}$ and $\textbf{T}_{t_2}$ are called transfer matrices. Because of such a nested relationship, the cluster bases only need to be stored for leaf clusters. For non-leaf clusters, only transfer matrices need to be stored. The ${\cal H}^2$-matrix is stored in a tree structure, with the size of leaf-level clusters denoted by $leafsize$. The number of blocks formed by a single cluster at each tree level is bounded by a constant $C_{sp}$.
In an ${\cal H}^2$-matrix, a large matrix block consisting of \textbf{F} and \textbf{R} is called a nonleaf block \textbf{NL}. As an example, a four-level block cluster ${\cal H}^2$-tree is illustrated in Fig. 1 (a), where the green link connects a row cluster with a column cluster, which form an admissible block, and the red links are for inadmissible blocks. The resultant ${\cal H}^2$-matrix is shown in Fig. 1 (b), where the admissible blocks are marked in green and the inadmissible blocks are marked in red. 

\section{Proposed ${\cal H}^2$ Matrix-Matrix Product Algorithm\textemdash Leaf Level}
\begin{figure}[t]
\centering
\includegraphics[width=0.48\textwidth]{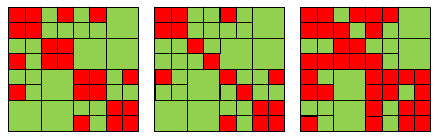}
\caption{An ${\cal H}^2$-matrix structure. (a)  $\textbf{A}_{{\cal H}^2}$. (b) $\textbf{B}_{{\cal H}^2}$. (c) $\textbf{C}_{{\cal H}^2}$.}
\end{figure}
To compute $\textbf{A}_{{\cal H}^2} \times \textbf{B}_{{\cal H}^2} =\textbf{C}_{{\cal H}^2}$, unlike the existing ${\cal H}^2$ formatted MMP \cite{BormH2MMP}, which is recursive, we propose to perform a one-way tree traversal from leaf level all the way up to the minimum level that has admissible blocks. Here, the tree is inverted with root level at level 0. While doing the multiplications at each level, we instantaneously compute the new row and column cluster bases based on prescribed accuracy to represent the product matrix accurately. We will use the ${\cal H}^2$-matrices shown in Fig. 2 to illustrate the proposed algorithm, but the algorithm is valid for any ${\cal H}^2$-matrix. The structures of $\textbf{A}_{{\cal H}^2}$, $\textbf{B}_{{\cal H}^2}$, and $\textbf{C}_{{\cal H}^2}$ matrices, i.e., which block is admissible and which is inadmissible, are determined based on the admissibility condition given in (\ref{admCond}). During the product calculation, we will keep the structure of product $\textbf{C}_{{\cal H}^2}$ matrix while achieving prescribed accuracy. %However, the algorithm can easily modified to another version where the the product structure will be changed due to the two multipliers as shown in \cite{MiaoAPSH2MMP}.
In this section, we detail proposed algorithm for leaf-level multiplications.

\begin{figure}[b]
\centering
\includegraphics[width=0.48\textwidth]{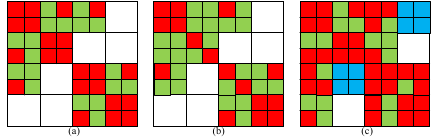}
\caption{${\cal H}^2$-matrix at leaf level. (a)  $\textbf{A}_{{\cal H}^2}^{L}$. (b) $\textbf{B}_{{\cal H}^2}^{L}$. (c) $\textbf{C}_{{\cal H}^2}^{L}$.}
\end{figure}
We start from leaf level $(l=L)$. Let \textbf{F} denote an inadmissible block, which is stored as a full matrix, and \textbf{R} be an admissible block. At leaf level, there are in total four matrix-matrix multiplication cases, i.e., 
\begin{itemize}
\item Case-1: $\textbf{F}^{\textbf{A}}$ $\times$ $\textbf{F}^{\textbf{B}}$ 
\item Case-2: $\textbf{F}^{\textbf{A}}$ $\times$ $\textbf{R}^{\textbf{B}}$
\item Case-3: $\textbf{R}^{\textbf{A}}$ $\times$ $\textbf{F}^{\textbf{B}}$
\item Case-4: $\textbf{R}^{\textbf{A}}$ $\times$ $\textbf{R}^{\textbf{B}}$
\end{itemize}
The resulting matrix block in $\textbf{C}$ is of two kinds: First, full matrix block, denoted by $\textbf{F}^{\textbf{C}}$, marked in red in Fig. 3(c); Second, admissible block of leaf size, which could be located at leaf level, denoted by $\textbf{R}^{\textbf{C},L}$ as marked in green in Fig. 3(c); which could also appear as a subblock in the non-leaf level $l$ as marked in blue in Fig. 3(c). The blue blocks in Fig. 3(c) are only for temporary storage, which will be changed to green admissible blocks during the upper level multiplication to preserve the structure of $\textbf{C}_{{\cal H}^2}$ matrix. The white blocks in Fig. 3 denote those blocks that are not involved in the leaf level multiplication. Next we show how to perform each matrix-matrix multiplication based on the two kinds of target blocks.
%denoted by $\textbf{R}^{\textbf{C},L}$. Notice that if the product block is inside an admissible block at nonleaf levels, that computation is performed at the corresponding non-leaf level of the target admissible block. 

\subsection{Product is an inadmissible block (full matrix) in $\textbf{C}$}
If the product matrix is a full block $\textbf{F}^{\textbf{C}}$, we can perform the four cases of multiplications exactly as they are by full matrix multiplications. For the admissible leaf blocks in four cases, we convert them into full matrices and then compute products. Since the size of these matrices is of $leafsize$, a user-defined constant, the computational cost is constant for each of such computations.

\subsection{Product is an admissible block in $\textbf{C}$}
\label{sec:leafAdm}
If the product is admissible in $\textbf{C}$ whether it is a leaf-level block or a subblock of a non-leaf admissible block, case-4 can be performed as it is since the product matrix is obviously admissible, which also preserves the original row and column cluster bases. In other words, the row cluster basis of $\textbf{A}$ is that of $\textbf{C}$; and the column cluster basis of $\textbf{B}$ is kept in $\textbf{C}$. To see this point clearly, we can write
\begin{equation}\label{case4}
\textbf{case-4: }
\textbf{R}_{i,j}^{\textbf{A}} \times \textbf{R}_{j,k}^{\textbf{B}}=\textbf{V}_{i_r}^{\textbf{A}}\textbf{S}_{i,j}^{\textbf{A}} (\textbf{V}_{j_c}^{\textbf{A}})^T \times \textbf{V}_{j_r}^{\textbf{B}}\textbf{S}_{j,k}^{\textbf{B}} (\textbf{V}_{k_c}^{\textbf{B}})^{T},
\end{equation}
where subscripts $i$, $j$, and $k$ denote cluster index, subscript $r$ denotes the corresponding cluster is a row cluster, whereas $c$ denotes the cluster is a column cluster. For example, $\textbf{V}_{i_r}^{\textbf{A}}$ denotes the cluster basis of row cluster $i$ in $\textbf{A}$, and $\textbf{V}_{k_c}^{\textbf{B}}$ denotes the cluster basis of column cluster $k$ in $\textbf{B}$. Eqn. (\ref{case4}) can be written in short as
\begin{equation}
\textbf{R}_{i,j}^{\textbf{A}} \times \textbf{R}_{j,k}^{\textbf{B}}=\textbf{V}_{i_r}^{\textbf{A}}\textbf{S}_{i,k}^{\textbf{C}}(\textbf{V}_{k_c}^{\textbf{B}})^{T},
\end{equation}
in which $\textbf{S}_{i,k}^{\textbf{C}}$ is the part in between the two cluster bases, which denotes the coupling matrix of the product admissible block in $\textbf{C}$. Clearly, this case of multiplication does not change the original row and column cluster bases.

For the other three cases, in existing MMP algorithms, a formatted multiplication is performed, which is done in the same way as case-4, i.e., using the original cluster bases of \textbf{A} and \textbf{B} or pre-assumed bases as the cluster bases of the product block. This obviously can be inaccurate since cases-1, 2, and 3, if performed as they are, would result in different cluster bases in the product matrix, which cannot be assumed. Specifically, case-1 results in a  different row as well as column cluster bases in the product admissible block because
\begin{equation}\label{case1}
\textbf{case-1: } \textbf{F}_{i,j}^{\textbf{A}} \times \textbf{F}_{j,k}^{\textbf{B}};
\end{equation}
case-2 yields a different row cluster basis since
\begin{equation}\label{case2}
\textbf{case-2: } \textbf{F}_{i,j}^{\textbf{A}} \times \textbf{R}_{j,k}^{\textbf{B}} = (\textbf{F}_{i,j}^{\textbf{A}}\textbf{V}_{j_r}^{\textbf{B}}) \times \textbf{S}_{j,k}^{\textbf{B}} \times (\textbf{V}_{k_c}^{\textbf{B}})^{T};
\end{equation}
whereas case-3 results in a different column cluster basis in the product admissible block, because
\begin{equation}\label{case3}
\textbf{case-3: } \textbf{R}_{i,j}^{\textbf{A}} \times \textbf{F}^{\textbf{B}}_{j,k} = \textbf{V}_{i_r}^{\textbf{A}} \times\textbf{S}_{i,j}^{\textbf{A}}\times \left((\textbf{V}_{j_c}^{\textbf{A}})^{T}\textbf{F}^{\textbf{B}}_{j,k}\right).
\end{equation}
If we do not update the cluster bases in the product matrix, the accuracy of the multiplication is not controllable. Therefore, in the proposed algorithm, we update row and column cluster bases for multiplication cases 1, 2, and 3 based on prescribed accuracy. We also have to do so with the nested property taken into consideration so that the computation at nonleaf levels can be performed efficiently.

For case-1, both row and column cluster bases of the product block need to be updated. For case-2, we need to use $\textbf{F}_{i,j}^{\textbf{A}}\textbf{V}_{j_r}^{\textbf{B}}$ to update the original row cluster basis $\textbf{V}_{i_r}^{\textbf{A}}$. For case-3, we need to use $(\textbf{V}_{j_c}^{\textbf{A}})^{T} \textbf{F}_{j,k}^{\textbf{B}}$ to update column cluster basis $\textbf{V}_{k_c}^{\textbf{B}}$. Since there are many case-1, 2 and 3 products encountered at the leaf level for the same row or column cluster, we develop the following algorithm to systematically update the cluster bases. In this procedure, \textbf{we also have to take the computation at all nonleaf levels into consideration so that the changed cluster bases at the leaf level can be reused at the nonleaf levels.} To achieve this goal, when we update the cluster basis due to the case-1, 2, and 3 multiplications associated with this cluster, not only we consider the product admissible block in the leaf level, but also the admissible blocks at all nonleaf levels. In other words, when computing $\textbf{A}_{i,j}$ multiplied by $\textbf{B}_{j,k}$, if the $\textbf{C}_{i,k}$ block is part of a non-leaf admissible block, we will take the corresponding multiplication into account to update the cluster bases. The detailed algorithms are as follows.

\subsection{Computation of new cluster bases in matrix product $\textbf{C}_{{\cal H}^2}$}

First, we show how to calculate the new row cluster bases of $\textbf{C}_{{\cal H}^2}$. Take an arbitrary row cluster $i$ as an example, let its cluster basis in $\textbf{C}$ be denoted by $\textbf{V}^{\textbf{C}}_{i_r}$. This cluster basis is affected by both case-1 and case-2 multiplications, as analyzed in the above. We first find all the case-1 multiplications associated with cluster $i$, i.e., all $\textbf{F}_{i,j}^{\textbf{A}} \times \textbf{F}_{j,k}^{\textbf{B}}$ whose product block $\textbf{C}_{i,k}$ is admissible. Again, notice that the $\textbf{C}_{i,k}$ can be either admissible at leaf level or be part of a non-leaf admissible block. For any cluster $i$, the number of $\textbf{F}_{i,j}^{\textbf{A}}$ is bounded by constant $C_{sp}$, since the number of inadmissible blocks that can be formed by a cluster is bounded by $C_{sp}$. For the same reason, the number of $\textbf{F}_{j,k}^{\textbf{B}}$ for cluster $j$ is also bounded by constant $C_{sp}$. Hence,  the total number of $\textbf{F}_{i,j}^{\textbf{A}} \times \textbf{F}_{j,k}^{\textbf{B}}$ multiplications is bounded by $C_{sp}^{2}$, thus also a constant. Then we calculate the Gram matrix sum of these products as:
\begin{equation}\label{Gr1L}
\textbf{G}_{i_{r1}}^{\textbf{C},L} =\sum_{j=1}^{O(C_{sp})}\sum_{k=1}^{O(C_{sp})}( \textbf{F}_{i,j}^{\textbf{A}} \textbf{F}_{j,k}^{\textbf{B}})(\textbf{F}_{i,j}^{\textbf{A}} \textbf{F}_{j,k}^{\textbf{B}})^{H},
\end{equation}
in which superscript $H$ denotes a Hermitian matrix. We also find all case-2 products associated with cluster $i$, which is the number of $\textbf{F}_{i,j}^{\textbf{A}}$ formed by cluster $i$ at leaf level in $\textbf{A}_{{\cal H}^2}$. This is also bounded by $C_{sp}$. Since in case-2 products, $\textbf{F}_{i,j}^{\textbf{A}}$ is multiplied by an admissible block in $\textbf{B}$, and hence $\textbf{V}_{j_r}^{\textbf{B}}$, we compute
\begin{equation}\label{Gr2L}
\textbf{G}_{i_{r2}}^{\textbf{C},L}=\sum_{j=1}^{O(C_{sp})}(\textbf{F}_{i,j}^{\textbf{A}}\textbf{V}_{j_r}^{\textbf{B}})(\textbf{F}_{i,j}^{\textbf{A}}\textbf{V}_{j_r}^{\textbf{B}})^{H},
\end{equation}
which incorporates all of the new cluster bases information due to case-2 products. 

For case-3 and case-4 multiplications, the row cluster bases of $\textbf{A}_{{\cal H}^2}$ matrix are kept to be those of $\textbf{C}$. So we account for the contribution of $\textbf{V}_{i_r}^{\textbf{A}}$ as
\begin{equation}\label{Gr3L}
\textbf{G}_{i_{r3}}^{\textbf{C},L}=\textbf{V}_{i_{r}}^{\textbf{A}}(\textbf{V}_{i_{r}}^{\textbf{A}})^{H}. 
\end{equation}
The column space spanning $\textbf{G}_{i_{r1}}^{\textbf{C},L}$, $\textbf{G}_{i_{r2}}^{\textbf{C},L}$ and $\textbf{G}_{i_{r3}}^{\textbf{C},L}$ would be the new cluster basis of $i$, since it takes both the original cluster basis and the change to the cluster basis due to matrix products into consideration. Since the magnitude of the three matrices may differ greatly, we normalize them before summing them up so that each component is captured. We thus obtain
\begin{equation}
\textbf{G}_{i_{r3}}^{\textbf{C},L}=\widehat{\textbf{G}_{i_{r1}}^{\textbf{C},L}}+\widehat{\textbf{G}_{i_{r2}}^{\textbf{C},L}}+\widehat{\textbf{G}_{i_{r3}}^{\textbf{C},L}}.
\end{equation}
The \ $\widehat{}$  \ above $\textbf{G}_{i_{r1}}^{\textbf{C},L}$, $\textbf{G}_{i_{r2}}^{\textbf{C},L}$ and $\textbf{G}_{i_{r3}}^{\textbf{C},L}$ denotes a normalized matrix. We then perform an SVD on $\textbf{G}_{i_{r3}}^{\textbf{C},L}$ to obtain the row cluster bases for cluster $i$ of $\textbf{C}_{{\cal H}^2}$ based on prescribed accuracy $\epsilon_{trunc}$. The singular vectors whose normalized singular values are greater than $\epsilon_{trunc}$ make the new row cluster basis $\textbf{V}_{i_r}^{\textbf{C}}$. It can be used to accurately represent the admissible blocks related to cluster $i$ in $\textbf{C}_{{\cal H}^2}$. Here, notice that the proposed algorithm for computing matrix-product cluster bases keeps nested property of $\textbf{V}_{i_r}^{\textbf{C}}$. This is because the Gram matrix sums in (\ref{Gr1L}), (\ref{Gr2L}) and (\ref{Gr3L}) take the upper level admissible products into account.

To compute the column cluster bases in $\textbf{C}_{{\cal H}^2}$, the steps are similar to the row cluster basis computation. We account for the contributions from all the four cases of products to compute column cluster bases. As can be seen from (\ref{case1}) and (\ref{case3}), in case-1 and case-3 products, the column cluster bases are changed from the original ones; whereas in case-2 and case-4 products, the column cluster bases are kept the same as those in $\textbf{B}$.

Consider an arbitrary column cluster $k$ in $\textbf{C}_{{\cal H}^2}$. We find all of the case-1 products associated with $k$, which is $\textbf{F}_{i,j}^{\textbf{A}} \times \textbf{F}_{j,k}^{\textbf{B}}$ with target $\textbf{C}_{i,k}$ being admissible either at the leaf or non-leaf level. The number of such multiplications is bounded by $C_{sp}^{2}$. We then compute the sum of their Gram matrices as:
\begin{equation}\label{case1-ch}
\textbf{G}_{k_{c1}}^{\textbf{C},L} =\sum_{i=1}^{O(C_{sp})} \sum_{j=1}^{O(C_{sp})}(\textbf{F}_{i,j}^{\textbf{A}} \textbf{F}_{j,k}^{\textbf{B}})^{T}(\textbf{F}_{i,j}^{\textbf{A}} \textbf{F}_{j,k}^{\textbf{B}})^{*}.
\end{equation}
Here, the superscript $^*$ denotes a complex conjugate. We also find all of the case-3 products associated with $k$, which is $\textbf{R}_{i,j}^{\textbf{A}} \times \textbf{F}_{j,k}^{\textbf{B}}$ with target $\textbf{C}_{i,k}$ being admissible either at the leaf or non-leaf level. Hence, the new column cluster basis takes a form of $(\textbf{V}_{j_c}^{\textbf{A}})^{T} \times \textbf{F}_{j,k}^{\textbf{B}}$. The number of such multiplications is also bounded by $C_{sp}$. The sum of their Gram matrices can be computed as:
\begin{equation}\label{case2-ch}
\textbf{G}_{k_{c2}}^{\textbf{C},L}=\sum_{j=1}^{O(C_{sp})}\left((\textbf{V}_{j_c}^{\textbf{A}})^{T}\textbf{F}_{j,k}^{\textbf{B}}\right)^{T}\left((\textbf{V}_{j_c}^{\textbf{A}})^{T}\textbf{F}_{j,k}^{\textbf{B}}\right)^{*}. 
\end{equation}
For case-2 and case-4 products, the original column cluster bases of $\textbf{B}_{{\cal H}^2}$ are kept in $\textbf{C}_{{\cal H}^2}$, hence, we compute 
\begin{equation}\label{case3-ch}
\textbf{G}_{k_{c3}}^{\textbf{C},L}=\textbf{V}_{k_{c}}^{\textbf{B}}(\textbf{V}_{k_{c}}^{\textbf{B}})^{H}. 
\end{equation}
We also normalize these three Gram matrices $\textbf{G}_{k_{c1}}^{\textbf{C},L}$, $\textbf{G}_{k_{c2}}^{\textbf{C},L}$ and $\textbf{G}_{k_{c3}}^{\textbf{C},L}$ and sum them up as:
\begin{equation}
\textbf{G}_{k_{c}}^{\textbf{C},L}=\widehat{\textbf{G}_{k_{c1}}^{\textbf{C},L}}+\widehat{\textbf{G}_{k_{c2}}^{\textbf{C},L}}+\widehat{\textbf{G}_{k_{c3}}^{\textbf{C},L}}.
\end{equation}
We then perform an SVD on this $\textbf{G}_{k_{c}}^{\textbf{C},L}$ and truncate the singular values based on prescribed accuracy $\epsilon_{trunc}$ to obtain the column cluster bases $\textbf{V}_{k_c}^{\textbf{C}}$ for cluster $k$. Now this new column cluster basis $\textbf{V}_{k_{c}}^{\textbf{C}}$ can be used to accurately represent the admissible blocks formed by column cluster $k$ in $\textbf{C}_{{\cal H}^2}$. 

\subsection{Computation of the four cases of multiplications with the product block being admissible}
After computing the new row and column cluster bases of the product matrix, for the multiplication whose target is an admissible block described in Section \ref{sec:leafAdm}, the computation becomes the coupling matrix computation since the cluster bases have been generated. For the four cases of multiplications, their coupling matrices have the following expressions:
\begin{equation}\label{leafcoup}
\textbf{S}_{i,k}^{\textbf{C}}= \begin{cases} (\textbf{V}_{i_r}^{\textbf{C}})^{H} \textbf{F}_{i,j}^{\textbf{A}} \textbf{F}_{j,k}^{\textbf{B}} (\textbf{V}_{k_c}^{\textbf{C}})^{*} &\text{case-1} \\
(\textbf{V}_{i_r}^{\textbf{C}})^{H} \textbf{F}_{i,j}^{\textbf{A}} \textbf{V}_{j_r}^{\textbf{B}}  \textbf{S}_{j,k}^{\textbf{B}}(\textbf{V}_{k_c}^{\textbf{B}})^{T}(\textbf{V}_{k_{c}}^{\textbf{C}})^{*} &\text{case-2} \\
(\textbf{V}_{i_{r}}^{\textbf{C}})^{H}\textbf{V}_{i_{r}}^{\textbf{A}}\textbf{S}_{i,j}^{\textbf{A}} (\textbf{V}_{j_c}^{\textbf{A}})^{T} \textbf{F}_{j,k}^{\textbf{B}} (\textbf{V}_{k_{c}}^{\textbf{C}})^{*} &\text{case-3} \\
(\textbf{V}_{i_{r}}^{\textbf{C}})^{H}\textbf{V}_{i_{r}}^{\textbf{A}}\textbf{S}_{i,j}^{\textbf{A}}\textbf{B}_{j}\textbf{S}_{j,k}^{\textbf{B}}(\textbf{V}_{k_{c}}^{\textbf{B}})^{T}(\textbf{V}_{k_{c}}^{\textbf{C}})^{*} &\text{case-4}.
\end{cases}
\end{equation}
The resulting admissible blocks in $\textbf{C}_{{\cal H}^2}$ are nothing but $\textbf{R}_{i,k}^{\textbf{C}}=\textbf{V}_{i_r}^{\textbf{C}} \times \textbf{S}_{i,k}^{\textbf{C}} \times (\textbf{V}_{k_c}^{\textbf{C}})^{T}$. 

In (\ref{leafcoup}), the $\textbf{B}_{j}$ is the cluster bases product, which is as shown below:
\begin{equation}
\textbf{B}_{j}=(\textbf{V}_{j_c}^{\textbf{A}})^{T} \times \textbf{V}_{j_r}^{\textbf{B}}.
\end{equation}
Since it is only related to the original cluster bases, it can be prepared in advance before the MMP computation. Using the nested property of the cluster bases, $\textbf{B}_{j}$ can be computed in linear time for all clusters $j$, be $j$ a leaf or a non-leaf cluster.

In (\ref{leafcoup}), the $(\textbf{V}_{i_r}^{\textbf{C}})^{H}\textbf{V}_{i_r}^{\textbf{A}}$ is simply the projection of the original row cluster basis of $\textbf{A}$ onto the new cluster basis of the product matrix $\textbf{C}$. Similarly, $ (\textbf{V}_{k_{c}}^{\textbf{B}})^{T}(\textbf{V}_{k_{c}}^{\textbf{C}})^{*}$ denotes the projection of the original column cluster basis of $\textbf{B}$  onto the newly generated column cluster basis in  $\textbf{C}$. The two cluster basis projections can also be computed for every leaf cluster after the new cluster bases have been generated. Hence, we compute
\begin{equation} \label{leafproj}
  \begin{array}{l}
    \textbf{P}_{i}^{\textbf{A}}=(\textbf{V}_{i_{r}}^{\textbf{C}})^{H}\textbf{V}_{i_{r}}^{\textbf{A}}; \\
   \textbf{P}_{k}^{\textbf{B}}= (\textbf{V}_{k_{c}}^{\textbf{B}})^{T}(\textbf{V}_{k_{c}}^{\textbf{C}})^{*}
  \end{array}
\end{equation}
for each leaf row cluster $i$, and each column leaf cluster $k$. In this way, it can be reused without recomputation for each admissible block formed by $i$ or $k$. 

In (\ref{leafcoup}), we can also see that the $\textbf{F}$ block is front and back multiplied by cluster bases. It can be viewed as an $\textbf{F}$ block collected based on the front (row) and back (column) cluster bases, which becomes a matrix of rank size. Specifically, in (\ref{leafcoup}), there are three kinds of collected blocks
\begin{equation}
  \begin{array}{l}\label{Fcollect}
    (\textbf{F}_{i,j}^{\textbf{A}}\textbf{F}_{j,k}^{\textbf{B}})_{coll.}=(\textbf{V}_{i_r}^{\textbf{C}})^{H}(\textbf{F}_{i,j}^{\textbf{A}} \textbf{F}_{j,k}^{\textbf{B}})\left(\textbf{V}_{k_{c}}^{\textbf{C}}\right)^* \\ 
    (\textbf{F}_{i,j}^{\textbf{A}})_{coll.}=(\textbf{V}_{i_r}^{\textbf{C}})^{H}\textbf{F}_{i,j}^{\textbf{A}}\textbf{V}_{j_r}^{\textbf{B}} \\ 
    (\textbf{F}_{j,k}^{\textbf{B}})_{coll.}=(\textbf{V}_{j_c}^{\textbf{A}})^{T}\textbf{F}_{j,k}^{\textbf{B}}\left(\textbf{V}_{k_c}^{\textbf{C}}\right)^{*},
  \end{array}
\end{equation}
which is used in case-1, 2, and 3 multiplication respectively. 

As can be seen from (\ref{leafcoup}), the case-1 multiplication with an admissible block being the target can be performed by first computing the full-matrix product, and then collecting the product onto the new row and column cluster bases of the product matrix. This collect operation is accurate because the newly generated row and column cluster bases have taken such a case-1 multiplication into consideration when being generated. %In contrast, in the algorithm of \cite{}, the collect operation is performed using the original cluster bases, and hence the accuracy cannot be controlled. 
 As for the case-2 multiplication, as can be seen from (\ref{leafcoup}), we can use the $\textbf{F}_{i,j}$ collected based on the new row cluster basis and the original column cluster basis, the size of which is rank, to multiply the coupling matrix of $\textbf{S}_{j,k}$, and then multiply the column basis projection matrix since the column bases have been changed. Similarly, for case-3, we use the collected block $ (\textbf{F}_{j,k}^{\textbf{B}})_{coll.}$, and front multiply it by the coupling matrix of $\textbf{S}_{i,j}$, and then front multiply a row cluster basis transformation matrix. As for case-4, we multiply the coupling matrix of $\textbf{A}$'s admissible block by the cluster basis product, and then by the coupling matrix of $\textbf{B}$'s admissible block. Since the row and column cluster bases have been changed to account for the other cases of multiplications, at the end, we need to front and back multiply the cluster basis transformation matrices to complete the computation of case-4. Summarizing the aforementioned, the coupling matrix in (\ref{leafcoup}) can be efficiently computed as
\begin{equation}\label{leafcoup_simple}
\textbf{S}_{i,k}^{\textbf{C}}= \begin{cases} (\textbf{V}_{i_r}^{\textbf{C}})^{H} \textbf{F}_{i,j}^{\textbf{A}} \textbf{F}_{j,k}^{\textbf{B}} (\textbf{V}_{k_c}^{\textbf{C}})^{*} &\text{case-1} \\
(\textbf{F}_{i,j}^{\textbf{A}})_{coll.}\textbf{S}_{j,k}^{\textbf{B}}\textbf{P}_{k}^{\textbf{B}} &\text{case-2} \\ 
\textbf{P}_{i}^{\textbf{A}}\textbf{S}_{i,j}^{\textbf{A}}(\textbf{F}_{j,k}^{\textbf{B}})_{coll.} &\text{case-3} \\
\textbf{P}_{i}^{\textbf{A}}\textbf{S}_{i,j}^{\textbf{A}}\textbf{B}_{j}\textbf{S}_{j,k}^{\textbf{B}}\textbf{P}_{k}^{\textbf{B}} &\text{case-4}.
\end{cases}
\end{equation}

\subsection{Summary of overall algorithm at leaf level}
Here, we conclude all the operations related to leaf level computation when the target is an admissible block:
\begin{enumerate}
\item Prepare cluster bases product $\textbf{B}$;
\item Compute all the leaf-level row and column cluster bases of product matrix $\textbf{C}_{{\cal H}^2}$;
\item Collect the \textbf{F} blocks in $\textbf{A}_{{\cal H}^2}$ and $\textbf{B}_{{\cal H}^2}$ based on the new row and/or column cluster bases, also prepare cluster bases transformation matrix $\textbf{P}$ ;
\item Perform four cases of multiplications.
\end{enumerate}

After leaf level multiplications, we need to merge four coupling matrices at a non-leaf level admissible block, as shown by the blue blocks in Fig. 3 (c). These matrices correspond to the multiplication case of a nonleaf block $\textbf{NL}$ multiplied by a nonleaf block $\textbf{NL}$ generating an admissible block at next level. The merged block is the coupling matrix of this next-level admissible block. It will be used to update next level transfer matrices. The details will be given in next section.

\section{Proposed ${\cal H}^2$ Matrix-Matrix Product Algorithm\textemdash Non-Leaf Level} \label{nonleafcomp}
\begin{figure}[b]
\centering
\includegraphics[width=0.45\textwidth]{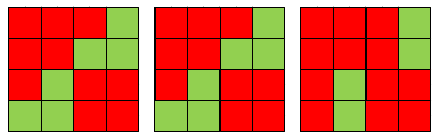}
\caption{${\cal H}^2$-matrix block at non-leaf level $(L-1)$. (a)  $\textbf{A}_{{\cal H}^2}^{L-1}$. (b) $\textbf{B}_{{\cal H}^2}^{L-1}$. (c) $\textbf{C}_{{\cal H}^2}^{L-1}$.}
\end{figure}
After finishing the leaf level multiplication, we proceed to non-leaf level multiplications. In Fig. 4, we use level $L-1$ as an example to illustrate $\textbf{A}_{{\cal H}^2}^{L-1}$, $\textbf{B}_{{\cal H}^2}^{L-1}$, and $\textbf{C}_{{\cal H}^2}^{L-1}$. 

At a nonleaf level $l$, there are also in total four matrix-matrix multiplication cases, i.e., 
\begin{itemize}
\item Case-1: $\textbf{NL}^{\textbf{A}}$ $\times$ $\textbf{NL}^{\textbf{B}}$
\item Case-2: $\textbf{NL}^{\textbf{A}}$ $\times$ $\textbf{R}^{\textbf{B}}$
\item Case-3: $\textbf{R}^{\textbf{A}}$ $\times$ $\textbf{NL}^{\textbf{B}}$
\item Case-4: $\textbf{R}^{\textbf{A}}$ $\times$ $\textbf{R}^{\textbf{B}}$,
\end{itemize}
where $\textbf{NL}$ denotes a non-leaf block. The resulting matrix block in $\textbf{C}$ is also of two kinds: 1) non-leaf block $\textbf{NL}$ at this level, marked in red in Fig. 4 (c), and 2) admissible block $\textbf{R}$, marked in green in Fig. 4 (c). Next we show how to perform each case of multiplications based on the two kinds of target blocks.

\subsection{Product is an $\textbf{NL}$ block in $\textbf{C}$}
\label{sec:nl-nl}
The $\textbf{NL}$ target block would not exist for a case-1 multiplication, since if a case-1 multiplication results in an $\textbf{NL}$ block, that computation should have been performed at previous level. As for the other three cases of multiplications, since at least one admissible block is present in the multipliers, the product must be an admissible block. Hence, we compute them as having an admissible block as the product, using the algorithm described in the following subsection, and associate the resulting admissible block with the $\textbf{NL}$ block. After the computation is done at all levels, we perform a backward split operation to split the admissible block associated with each $\textbf{NL}$ block to each leaf block of $\textbf{C}$ based on its structure.

\subsection{Product is an admissible block in $\textbf{C}$}
Similar to the leaf level, if the product is an admissible block in $\textbf{C}$ whether at the same non-leaf level or at an upper level, case-4 can be performed as it is since the product matrix is obviously admissible, which also preserves the original row and column cluster bases. We can write
\begin{align}\label{case4-nl}
&\textbf{case-4: }\nonumber\\
&\textbf{R}_{i,j}^{\textbf{A}} \times \textbf{R}_{j,k}^{\textbf{B}}=\nonumber\\
& \textbf{V}_{i^{ch}_r}^{\textbf{A}}\textbf{T}_{i_r}^{\textbf{A}}\textbf{S}_{i,j}^{\textbf{A}} (\textbf{T}_{j_c}^{\textbf{A}})^T(\textbf{V}_{j^{ch}_c}^{\textbf{A}})^T \times \textbf{V}_{j^{ch}_r}^{\textbf{B}}\textbf{T}_{j_r}^{\textbf{B}}\textbf{S}_{j,k}^{\textbf{B}} (\textbf{T}_{k_c}^{\textbf{B}})^T(\textbf{V}_{k^{ch}_c}^{\textbf{B}})^{T},
\end{align}
where $\textbf{T}$ denotes a transfer matrix, and superscript $ch$ denotes the two children clusters of the non-leaf cluster $i$. If the cluster bases at leaf level and the transfer matrices at non-leaf levels are kept the same as before, then the computation of (\ref{case4-nl}) is to calculate the coupling matrix at level $l$, which is
\begin{equation}
\textbf{S}_{i,k}^{\textbf{C}}=\textbf{S}_{i,j}^{\textbf{A}} (\textbf{B}_{j})\textbf{S}_{j,k}^{\textbf{B}}.
\end{equation}
It is a product of three small matrices whose size is the rank at this tree level. Rewriting (\ref{case4-nl}) as
\begin{align}\label{case4-2}
&\textbf{case-4: }\nonumber\\
&\textbf{R}_{i,j}^{\textbf{A}} \times \textbf{R}_{j,k}^{\textbf{B}}=\textbf{V}_{i^{ch}_r}^{\textbf{A}}\textbf{T}_{i_r}^{\textbf{A}}\textbf{S}_{i,k}^{\textbf{C}}(\textbf{T}_{k_c}^{\textbf{B}})^T(\textbf{V}_{k^{ch}_c}^{\textbf{B}})^{T}.
\end{align}
If we exclude the children cluster bases in the front and at the back, we can see that $\textbf{T}$ serves as the new cluster basis at this level. In other words, at a non-leaf level $l$, if we treat this level as the bottom level of the remaining tree, then the transfer matrix of the non-leaf cluster is nothing but the leaf cluster basis of the shortened tree.

Similar to the leaf-level computation, the other three cases of multiplications will result in a change of cluster basis in the matrix product. Specifically, case-1 results in a different row as well as column cluster bases in the product admissible block because
\begin{equation}\label{case1}
\textbf{case-1: } \textbf{NL}_{i,j}^{\textbf{A}} \times \textbf{NL}_{j,k}^{\textbf{B}};
\end{equation}
case-2 yields a different row cluster basis since
\begin{equation}\label{case2}
\textbf{case-2: } \textbf{NL}_{i,j}^{\textbf{A}} \times \textbf{R}_{j,k}^{\textbf{B}} = (\textbf{NL}_{i,j}^{\textbf{A}}\textbf{V}_{j_r}^{\textbf{B}}) \times \textbf{S}_{j,k}^{\textbf{B}} \times (\textbf{V}_{k_c}^{\textbf{B}})^{T};
\end{equation}
whereas case-3 results in a different column cluster basis in the product admissible block, because
\begin{equation}\label{case3}
\textbf{case-3: } \textbf{R}_{i,j}^{\textbf{A}} \times \textbf{NL}^{\textbf{B}}_{j,k} = \textbf{V}_{i_r}^{\textbf{A}} \times\textbf{S}_{i,j}^{\textbf{A}}\times \left((\textbf{V}_{j_c}^{\textbf{A}})^{T}\textbf{NL}^{\textbf{B}}_{j,k}\right).
\end{equation}
If we do not update the cluster bases in the product matrix, the accuracy of the multiplication is not controllable. However, if we update the cluster basis as they are, it is computationally very expensive since the matrix block size keeps increasing when we proceed from leaf level towards the root level. In addition to the cost of changing cluster bases, if we have to carry out the multiplications at each non-leaf level using the actual matrix block size, then the computation is also prohibitive. Therefore, the fast algorithm we develop here is to perform all computations using the rank size at each tree level, and meanwhile control the accuracy.

%Therefore, in the proposed algorithm, we update row and column cluster bases for multiplication cases 1, 2, and 3 based on prescribed accuracy. We also have to do so with the nested property taken into consideration so that the computation at nonleaf levels can be performed efficiently.

In the proposed algorithm, to account for the updates to the original matrix during the MMP procedure, the cluster bases of $\textbf{C}$ are computed level by level, which are manifested by the changed leaf cluster bases and the transfer matrices at nonleaf levels. At a non-leaf level, its children-level cluster bases have already been computed, and they are different from the original ones in $\textbf{A}$ and $\textbf{B}$. However, the new cluster bases have taken the upper-level multiplications into consideration. Hence, we can accurately represent the multiplication at the current non-leaf level using newly generated children cluster bases.

Take case-1 product as an example, where we perform $\textbf{NL}_{i,j}^{\textbf{A}} \times \textbf{NL}_{j,k}^{\textbf{B}}$ obtaining an admissible $\textbf{R}_{i,k}^{\textbf{C}}$. %This is computed first as an $\textbf{NL}$ block, which is then merged as an admissible block $\textbf{R}$. Specifically, we split the $\textbf{R}_{i,k}^{\textbf{C}}$ into four admissible blocks, and compute each block at pevious level. The coupling matrices of the four admissible blocks are then merged together to become the  coupling matrix of enire $\textbf{R}_{i,k}^{\textbf{C}}$. This result obviously require a change of the transfer matrix of non-leaf cluster $i$ since the matrix product may not be represented by the original transfer matrix. 
We can accurately represent this product using the children cluster bases of $i$ and $k$ as follows:
\begin{align}\label{case1_nl}
&\textbf{case-1: }\nonumber\\
&\textbf{NL}_{i,j}^{\textbf{A}} \times \textbf{NL}_{j,k}^{\textbf{B}}=\nonumber\\
&\begin{bmatrix}
\textbf{V}_{i1_r}^{\textbf{C}} &  \\
 & \textbf{V}_{i2_r}^{\textbf{C}} \\
\end{bmatrix}(\textbf{NL}_{i,j}^{\textbf{A}}\textbf{NL}_{j,k}^{\textbf{B}})_{coll.} \begin{bmatrix}
(\textbf{V}_{k1_c}^{\textbf{C}})^T &  \\
 & (\textbf{V}_{k2_c}^{\textbf{C}})^T \\\end{bmatrix},
\end{align}
in which
\begin{align}\label{NL_col}
&(\textbf{NL}_{i,j}^{\textbf{A}}\textbf{NL}_{j,k}^{\textbf{B}})_{coll.}=\nonumber\\
&\begin{bmatrix}
(\textbf{V}_{i1_r}^{\textbf{C}})^H &  \\
 & (\textbf{V}_{i2_r}^{\textbf{C}})^H \\
\end{bmatrix}\left(\textbf{NL}_{i,j}^{\textbf{A}}\textbf{NL}_{j,k}^{\textbf{B}}\right)\begin{bmatrix}
(\textbf{V}_{k1_c}^{\textbf{C}})^* &  \\
 & (\textbf{V}_{k2_c}^{\textbf{C}})^* \\
\end{bmatrix}.
\end{align}
This collected block, $(\textbf{NL}_{i,j}^{\textbf{A}}\textbf{NL}_{j,k}^{\textbf{B}})_{coll.}$, is actually the coupling matrix merged from the four small coupling matrices computed at previous level, when dealing with the multiplication case of having a target block as a subblock in the upper-level admissible block. It can be written as
\begin{equation} 
(\textbf{NL}_{i,j}^{\textbf{A}}\textbf{NL}_{j,k}^{\textbf{B}})_{coll.}=
\begin{bmatrix} \label{nonleafA}
\textbf{S}_{i_1,k_1}^{\textbf{C}} &  \textbf{S}_{i_1,k_2}^{\textbf{C}} \\
\textbf{S}_{i_2,k_1}^{\textbf{C}} & \textbf{S}_{i_2,k_2}^{\textbf{C}}\\
\end{bmatrix}.
\end{equation}
Each of the four coupling matrices has been obtained at previous level. From (\ref{NL_col}), it is clear that using the nested property of the cluster bases, the collect operation does not need to start from leaf level, but using the four blocks obtained at previous one level. 

For case-2 product, it can also be accurately expanded in the space of the children row cluster bases, and hence
\begin{align}\label{case2_nl}
&\textbf{case-2: }\textbf{NL}_{i,j}^{\textbf{A}} \times \textbf{R}_{j,k}^{\textbf{B}}=\nonumber\\
&\begin{bmatrix}
\textbf{V}_{i1_r}^{\textbf{C}} &  \\
 & \textbf{V}_{i2_r}^{\textbf{C}} \\
\end{bmatrix}{\textbf{NL}_{i,j}^{\textbf{A}}}_{coll.}\textbf{T}_{j_r}^{\textbf{B}}\textbf{S}_{j,k}^{\textbf{B}}(\textbf{V}_{k_c}^{\textbf{B}})^T,
\end{align}
where
\begin{align}\label{case2_coll}
{\textbf{NL}_{i,j}^{\textbf{A}}}_{coll.}=\begin{bmatrix}
(\textbf{V}_{i1_r}^{\textbf{C}})^H &  \\
 & (\textbf{V}_{i2_r}^{\textbf{C}})^H \\
\end{bmatrix}\textbf{NL}_{i,j}^{\textbf{A}}\begin{bmatrix}
\textbf{V}_{j1_r}^{\textbf{B}} &  \\
 & \textbf{V}_{j2_r}^{\textbf{B}} \\
\end{bmatrix},
\end{align}
which is ${\textbf{NL}_{i,j}^{\textbf{A}}}$ collected based on the children's new row cluster bases in $\textbf{C}$ and the original column cluster bases in $\textbf{B}$. From (\ref{case2_nl}), it can be seen that if excluding the children cluster bases, then ${\textbf{NL}_{i,j}^{\textbf{A}}}_{coll.}\textbf{T}_{j_r}^{\textbf{B}}$ resembles the $\textbf{F}_{i,j}^{\textbf{A}}\textbf{V}_{j_r}^{\textbf{B}}$ in the leaf level case-2 product. In other words, if we treat the current non-leaf level as the leaf level, then ${\textbf{NL}_{i,j}^{\textbf{A}}}_{coll.}$ is equivalent to a full matrix block, whereas $\textbf{T}$ is the leaf cluster basis. An example of ${\textbf{NL}_{i,j}^{\textbf{A}}}_{coll.}$ block at level $(L-1)$ in $\textbf{A}_{{\cal H}^2}$ can be seen below:
\begin{equation} 
(\textbf{NL}_{i,j}^{\textbf{A}})_{coll.}=
\begin{bmatrix} \label{nonleafA}
(\textbf{F}_{i_1,j_1}^{\textbf{A}})_{coll.} &  \textbf{P}_{i_1}^{\textbf{A}}\textbf{S}_{i_1,j_2}^{\textbf{A}}\textbf{B}_{j_2} \\
\textbf{P}_{i_2}^{\textbf{A}}\textbf{S}_{i_2,j_1}^{\textbf{A}}\textbf{B}_{j_1} & (\textbf{F}_{i_2,j_2}^{\textbf{A}})_{coll.} \\
\end{bmatrix},
\end{equation}
which consists of collected full matrices whose expressions are shown in (\ref{Fcollect}), and projected coupling matrices of admissible blocks. Again, using the nested property of both new and original cluster bases, the collect operation does not need to start from leaf level, but using the four blocks obtained at previous one level. Each collect operation only costs $O(k_l)^3$, where $k_l$ is the rank at level $l$.

Since the cluster bases at the previous level have been computed, for case-1 and case-2 products at a non-leaf level, we only need to compute the center block associated with the current non-leaf level, and this computation can be carried out in the same way as how we carry out leaf-level computation, if we treat the current non-leaf level as the leaf level of the remaining tree. The same is true to case-3 product, where we have 
\begin{align}\label{case3_nl}
&\textbf{case-3: }\textbf{R}_{i,j}^{\textbf{A}} \times \textbf{NL}_{j,k}^{\textbf{B}}=\nonumber\\
&\textbf{V}_{i_r}^{\textbf{A}} \textbf{S}_{i,j}^{\textbf{A}}(\textbf{T}_{j_c}^{\textbf{A}})^T{\textbf{NL}_{j,k}^{\textbf{B}}}_{coll.}
\begin{bmatrix}
(\textbf{V}_{k1_c}^{\textbf{C}})^T &  \\
 & (\textbf{V}_{k2_c}^{\textbf{C}})^T \\
\end{bmatrix},
\end{align}
in which
\begin{align}\label{case3_coll}
{\textbf{NL}_{j,k}^{\textbf{B}}}_{coll.}=\begin{bmatrix}
(\textbf{V}_{j1_r}^{\textbf{A}})^T &  \\
 & (\textbf{V}_{j2_r}^{\textbf{A}})^T \\
\end{bmatrix}\textbf{NL}_{j,k}^{\textbf{B}}\begin{bmatrix}
(\textbf{V}_{k1_c}^{\textbf{C}})^* &  \\
 & (\textbf{V}_{k2_c}^{\textbf{C}})^* \\
\end{bmatrix}.
\end{align}
We can see that $(\textbf{T}_{j_c}^{\textbf{A}})^T{\textbf{NL}_{j,k}^{\textbf{B}}}_{coll.}$ resembles the $(\textbf{V}_{j_c}^{\textbf{A}})^T\textbf{F}_{j,k}^{\textbf{B}}$ in the leaf level case-3 product. An example of collected \textbf{NL} block in $\textbf{B}_{{\cal H}^2}$ is given as follows
\begin{equation}
(\textbf{NL}_{i,j}^{\textbf{B}})_{coll.}=
\begin{bmatrix}\label{nonleafB}
(\textbf{F}_{i_1,j_1}^{\textbf{B}})_{coll.}&  \textbf{B}_{i_1}\textbf{S}_{i_1,j_2}^{\textbf{B}}\textbf{P}_{j_2}^{\textbf{B}} \\
 \textbf{B}_{i_2}\textbf{S}_{i_2,j_1}^{\textbf{B}}\textbf{P}_{j_1}^{\textbf{B}} & (\textbf{F}_{i_2,j_2}^{\textbf{B}})_{coll.} \\
\end{bmatrix},
\end{equation}
which consists of collected full matrices whose expressions are shown in (\ref{Fcollect}), and projected coupling matrices of admissible blocks.

Since the cluster bases have been changed at previous level, we also represent the case-4 product using the new children cluster bases of $i$ and $k$, thus
\begin{align}\label{case4_nl}
&\textbf{case-4: }\nonumber\\
&\textbf{R}_{i,j}^{\textbf{A}} \times \textbf{R}_{j,k}^{\textbf{B}}=\nonumber\\
&\begin{bmatrix}
(\textbf{V}_{i1_r}^{\textbf{C}}) &  \\
 & (\textbf{V}_{i2_r}^{\textbf{C}}) \\
\end{bmatrix}\textbf{R}_{i,k,proj}^{\textbf{C}} \begin{bmatrix}
(\textbf{V}_{k1_c}^{\textbf{C}})^T &  \\
 & (\textbf{V}_{k2_c}^{\textbf{C}})^T \\\end{bmatrix},
\end{align}
and 
\begin{align}\label{case4_new}
&\textbf{R}_{i,k,proj}^{\textbf{C}}=\begin{bmatrix}
(\textbf{P}_{i1}^{\textbf{A}}) &  \\
 & (\textbf{P}_{i2}^{\textbf{A}}) \\
\end{bmatrix}(\textbf{T}_{i_r}^{\textbf{A}}\textbf{S}_{i,k}^{\textbf{C}}(\textbf{T}_{k_c}^{\textbf{B}})^T)\begin{bmatrix}
(\textbf{P}_{k1}^{\textbf{B}}) &  \\
 & (\textbf{P}_{k2}^{\textbf{B}}) \\
\end{bmatrix},
\end{align}
which can be written in short as
\begin{align}\label{case4_new}
&\textbf{R}_{i,k,proj}^{\textbf{C}}=\textbf{P}_{i^{ch}}^{\textbf{A}}(\textbf{T}_{i_r}^{\textbf{A}}\textbf{S}_{i,k}^{\textbf{C}}(\textbf{T}_{k_c}^{\textbf{B}})^T)\textbf{P}_{k^{ch}}^{\textbf{B}},
\end{align}
where $ch$ denotes children. Here, there is a cluster basis transformation matrix in the front and at the back.

%if we split the target admissible block into four small admissible blocks, we will see all the multiplications belong to the case of $\textbf{A}_{i_\xi,j_\xi} \times \textbf{B}_{j_\xi,k_\xi}$ generating an admissible $\textbf{C}_{i_\xi,k_\xi}$, where $i_\xi$, $j_\xi$, and $k_\xi$ with $\xi=1,2$ are the two children clusters of the nonleaf cluster $i$, $j$, and $k$. These multiplications have been used to generate the new row cluster basis of cluster $i_\xi$, and column cluster basis of $k_\xi$ at previous level. 
 
\subsection{Computation of the new non-leaf level transfer matrices in $\textbf{C}$}
If the target block is an admissible block at a nonleaf level, we need to represent it as $\textbf{R}_{t,s}=\textbf{T}_t\textbf{S}_{t,s}(\textbf{T}_s)^{T}$ in controlled accuracy. Hence, we need to calculate new row and column transfer matrices $\textbf{T}$ of product matrix $\textbf{C}_{{\cal H}^2}$. First, we introduce how to calculate the row transfer matrices. Similar to leaf level, case-1 and 2 products result in a change in the row cluster basis and hence row transfer matrix. Case-3 and 4 products do not require a change of transfer matrix if the cluster bases have not been changed at previous level. However, since the cluster bases have been changed at previous level, the transfer matrix requires an update as well. %But the update is nothing but multiplying the original transfer matrix by a cluster basis transformation matrix. %which also consists of three steps as similar to leaf level cluster bases calculation. The first step is using the admissible products at leaf level $L$ computation that located in nonleaf admissible blocks, the blue blocks in Fig. 3 (c). During leaf level computation, these blue blocks are stored as leaf level coupling matrix. 

For an arbitrary non-leaf cluster $i$, we first find all of the case-1 products associated with $i$. Each of such a product leads to a coupling matrix merged from the four coupling matrices obtained at previous level computation, denoted by $(\textbf{NL}_{i,j}^{\textbf{A}}\textbf{NL}_{j,k}^{\textbf{B}})_{coll.}$. Using them, we calculate the Gram matrix sum as:
%% The number of NL^c is not bounded by Csp??????
\begin{equation}\label{case1-ch-nl}
\textbf{G}_{i_{r1}}^{\textbf{C},l}=\sum_{\#(i,k)=1}^{O(C_{sp}^2)}(\textbf{NL}_{i,j}^{\textbf{A}}\textbf{NL}_{j,k}^{\textbf{B}})_{coll.}((\textbf{NL}_{i,j}^{\textbf{A}}\textbf{NL}_{j,k}^{\textbf{B}})_{coll.})^{H}.
\end{equation}

The second step is to take case-2 multiplications at a non-leaf level into consideration for row transfer matrix calculation of product matrix $\textbf{C}_{{\cal H}^2}$. We find all the collected nonleaf blocks ${\textbf{NL}_{i,j}^{\textbf{A}}}_{coll.}$ of cluster $i$ at level $l$ in $\textbf{A}_{{\cal H}^2}$ matrix and multiply them with corresponding transfer matrices $\textbf{T}_{j_r}^{\textbf{B}}$ from $\textbf{B}_{{\cal H}^2}$ matrix. And we calculate the Gram matrix sum as
\begin{equation}\label{case2-ch-nl}
\textbf{G}_{i_{r2}}^{\textbf{C},l}=\sum_{j=1}^{O(C_{sp})}((\textbf{NL}_{i,j}^{\textbf{A}})_{coll.}\textbf{T}_{j_r}^{\textbf{B}})((\textbf{NL}_{i,j}^{\textbf{A}})_{coll.}\textbf{T}_{j_r}^{\textbf{B}})^{H}.
\end{equation}
Finally, we count the contributions from case-3 and case-4 products by computing
\begin{equation}\label{case3-ch-nl}
\textbf{G}_{i_{r3}}^{\textbf{C},l}=\textbf{P}_{i^{ch}}^{\textbf{A}}\textbf{T}_{i_{r}}^{\textbf{A}}(\textbf{T}_{i_{r}}^{\textbf{A}})^{H}(\textbf{P}_{i^{ch}}^{\textbf{A}})^H. 
\end{equation}
Again, we normalize these three Gram matrices and obtain
\begin{equation}
\textbf{G}_{i_{r}}^{\textbf{C},l}=\widehat{\textbf{G}_{i_{r1}}^{\textbf{C},l}}+\widehat{\textbf{G}_{i_{r2}}^{\textbf{C},l}} +\widehat{\textbf{G}_{i_{r3}}^{\textbf{C},l}}.
\end{equation}
We then calculate an SVD of this $\textbf{G}_{i_{r}}^{\textbf{C},l}$ and truncate the singular values based on prescribed accuracy $\epsilon_{trunc}$ to obtain row transfer matrix $\textbf{T}_{i_{r}}^{\textbf{C}}$ for cluster $i$ at nonleaf level. 

Similarly, we can compute the new column transfer matrices for non-leaf cluster $k$, which is $\textbf{T}_{k_{c}}^{\textbf{C}}$. The first part is 
\begin{equation}
\textbf{G}_{k_{c1}}^{\textbf{C},l} =\sum_{i=1}^{O(C_{sp})}((\textbf{NL}_{i,j}^{\textbf{A}}\textbf{NL}_{j,k}^{\textbf{B}})_{coll.})^{T}((\textbf{NL}_{i,j}^{\textbf{A}}\textbf{NL}_{j,k}^{\textbf{B}})_{coll.})^{*}.
\end{equation}
The second part is
\begin{equation}
\textbf{G}_{k_{c2}}^{\textbf{C},l}=\sum_{j=1}^{O(C_{sp})}((\textbf{T}_{j_c}^{\textbf{A}})^{T}(\textbf{NL}_{j,k}^{\textbf{B}})_{coll.})^{T}((\textbf{T}_{j_c}^{\textbf{A}})^{T}(\textbf{NL}_{j,k}^{\textbf{B}})_{coll.})^{*}. 
\end{equation}
The third part is 
\begin{equation}
\textbf{G}_{k_{c3}}^{\textbf{C},l}=(\textbf{P}_{k^{ch}}^{\textbf{B}})^T(\textbf{T}_{k_{c}}^{\textbf{B}})^{*}(\textbf{T}_{k_{c}}^{\textbf{B}})^{T}(\textbf{P}_{k^{ch}}^{\textbf{B}})^*. 
\end{equation}
Then we normalize the three Gram matrices and sum them up as 
\begin{equation}
\textbf{G}_{k_{c}}^{\textbf{C},l}=\widehat{\textbf{G}_{k_{c1}}^{\textbf{C},l}}+\widehat{\textbf{G}_{k_{c2}}^{\textbf{C},l}}+\widehat{\textbf{G}_{k_{c3}}^{\textbf{C},l}}.
\end{equation}
After we perform an SVD on $\textbf{G}_{k_{c}}^{\textbf{C},l}$ matrix and truncate  the singular values based on prescribed accuracy $\epsilon_{trunc}$, we get new column transfer matrix $\textbf{T}_{k_{c}}^{\textbf{C}}$. 

\subsection{Computation of the four cases of multiplications with the product block being admissible}

Now we obtain both row and column transfer matrices for product matrix $\textbf{C}_{{\cal H}^2}$, hence, the four multiplications become the computation of the coupling matrices, so that the admissible block at the current level has a form of $\textbf{R}_{t,s}=\textbf{T}_t\textbf{S}_{t,s}(\textbf{T}_s)^{T}$. The coupling matrix $\textbf{S}$'s calculation is similar to that of leaf level in (\ref{leafcoup}), which has the following expressions:
\begin{equation}\label{nonleafcoup}
\textbf{S}_{i,k}^{\textbf{C}}= \begin{cases} (\textbf{T}_{i_{r}}^{\textbf{C}})^{H}(\textbf{NL}_{i,j}^{\textbf{A}}\textbf{NL}_{j,k}^{\textbf{B}})_{coll.}(\textbf{T}_{k_{c}}^{\textbf{C}})^{*} &\text{case-1} \\
(\textbf{T}_{i_{r}}^{\textbf{C}})^{H} (\textbf{NL}_{i,j}^{\textbf{A}})_{coll.} \textbf{T}_{j_r}^{\textbf{B}}  \textbf{S}_{j,k}^{\textbf{B}}(\textbf{V}_{k_{c}}^{\textbf{B}})^{T}(\textbf{V}_{k_{c}}^{\textbf{C}})^{*} &\text{case-2} \\
(\textbf{V}_{i_{r}}^{\textbf{C}})^{H}\textbf{V}_{i_{r}}^{\textbf{A}}\textbf{S}_{i,j}^{\textbf{A}} (\textbf{T}_{j_c}^{\textbf{A}})^{T} (\textbf{NL}_{j,k}^{\textbf{B}})_{coll.} (\textbf{T}_{k_{c}}^{\textbf{C}})^{*} &\text{case-3} \\
(\textbf{V}_{i_{r}}^{\textbf{C}})^{H}\textbf{V}_{i_{r}}^{\textbf{A}}\textbf{S}_{i,j}^{\textbf{A}}\textbf{B}_{j}\textbf{S}_{j,k}^{\textbf{B}}(\textbf{V}_{k_{c}}^{\textbf{B}})^{T}(\textbf{V}_{k_{c}}^{\textbf{C}})^{*} &\text{case-4}.
\end{cases}
\end{equation}

Again, we should prepare some matrix products in advance so that we can achieve linear complexity MMP for constant rank ${\cal H}^2$-matrix. For nonleaf levels, the cluster bases product $\textbf{B}_{j}$ can be readily calculated using children's cluster bases based on the nested property. For example, given a nonleaf cluster $j$, we can generate $\textbf{B}_{j}$ by using the cluster bases product of its children clusters $j_1$ and $j_2$, which is shown as:
\begin{equation}\label{nonleafproduct}
\textbf{B}_{j}=(\textbf{T}_{{j_1}_c}^{\textbf{A}})^{T} \textbf{B}_{j_1} \textbf{T}_{{j_1}_r}^{\textbf{B}} +(\textbf{T}_{{j_2}_c}^{\textbf{A}})^{T} \textbf{B}_{j_2} \textbf{T}_{{j_2}_r}^{\textbf{B}}.
\end{equation}
Besides, since the cluster bases product $\textbf{B}_{j}$ only involve original cluster bases in $\textbf{A}_{{\cal H}^2}$ and $\textbf{B}_{{\cal H}^2}$ matrices, we can prepare the above $\textbf{B}_{j}$ for all leaf and nonleaf clusters before MMP algorithm. In addition, the nonleaf level cluster bases projection (transformation) can also be calculated using children's ones as shown in  (\ref{leafproj}). The formulas are given below: 
\begin{equation}\label{nonleafproj}
  \begin{aligned}
    \textbf{P}_{i}^{\textbf{A}} &=(\textbf{V}_{i_{r}}^{\textbf{C}})^{H}\textbf{V}_{i_{r}}^{\textbf{A}} \\
    &=(\textbf{T}_{{i_1}_r}^{\textbf{C}})^{H} \textbf{P}_{i_1}^{\textbf{A}} \textbf{T}_{{i_1}_r}^{\textbf{A}} +(\textbf{T}_{{i_2}_r}^{\textbf{C}})^{H} \textbf{P}_{i_2}^{\textbf{A}} \textbf{T}_{{i_2}_r}^{\textbf{A}}; \\
   \textbf{P}_{k}^{\textbf{B}}&=(\textbf{V}_{k_{c}}^{\textbf{B}})^{T}(\textbf{V}_{k_{c}}^{\textbf{C}})^{*} \\
   &= (\textbf{T}_{{k_1}_c}^{\textbf{B}})^{T} \textbf{P}_{k_1}^{\textbf{B}}(\textbf{T}_{{k_1}_c}^{\textbf{C}})^{*}+(\textbf{T}_{{k_2}_c}^{\textbf{B}})^{T} \textbf{P}_{k_2}^{\textbf{B}}(\textbf{T}_{{k_2}_c}^{\textbf{C}})^{*}.
  \end{aligned}
\end{equation}
We also compute the collected \textbf{NL} matrix block in $\textbf{A}_{{\cal H}^2}$ and $\textbf{B}_{{\cal H}^2}$ at current level $l$ by the following equation,
\begin{equation}
  \begin{array}{l}\label{NLcollect}
    (\textbf{NL}_{i,j}^{\textbf{A}})_{coll.}^{(l)}=(\textbf{T}_{i_{r}}^{\textbf{C}})^{H} (\textbf{NL}_{i,j}^{\textbf{A}})_{coll.}^{(l+1)} \textbf{T}_{j_r}^{\textbf{B}} \\ 
    (\textbf{NL}_{j,k}^{\textbf{B}})_{coll.}^{(l)}=(\textbf{T}_{j_c}^{\textbf{A}})^{T}(\textbf{NL}_{j,k}^{\textbf{B}})_{coll.}^{(l+1)}(\textbf{T}_{k_{c}}^{\textbf{C}})^{*}
  \end{array}
\end{equation}
where superscript $l$ denotes tree level. After we prepare the matrix products in (\ref{nonleafproduct}), (\ref{nonleafproj}), and (\ref{NLcollect}), we can proceed to calculate the coupling matrices in (\ref{nonleafcoup}) efficiently as:
\begin{equation}\label{nonleafcoup_simple}
\textbf{S}_{i,k}^{\textbf{C}}= \begin{cases} (\textbf{T}_{i_{r}}^{\textbf{C}})^{H}(\textbf{NL}_{i,j}^{\textbf{A}}\textbf{NL}_{j,k}^{\textbf{B}})_{coll.}^{(l)}(\textbf{T}_{k_{c}}^{\textbf{C}})^{*} &\text{case-1} \\
(\textbf{NL}_{i,j}^{\textbf{A}})_{coll.}^{(l)}\textbf{S}_{j,k}^{\textbf{B}}\textbf{P}_{k}^{\textbf{B}} &\text{case-2} \\ 
\textbf{P}_{i}^{\textbf{A}}\textbf{S}_{i,j}^{\textbf{A}}(\textbf{NL}_{j,k}^{\textbf{B}})_{coll.}^{(l)} &\text{case-3} \\
\textbf{P}_{i}^{\textbf{A}}\textbf{S}_{i,j}^{\textbf{A}}\textbf{B}_{j}\textbf{S}_{j,k}^{\textbf{B}}\textbf{P}_{k}^{\textbf{B}} &\text{case-4}.
\end{cases}
\end{equation}
All the coupling matrices calculation are performed in rank size $k_l$. So the computational cost is $O(k_l^3)$. After coupling matrices calculation in (\ref{nonleafcoup_simple}), all the admissible products at this nonleaf level multiplication can be represented as  $\textbf{R}_{i,j}^{\textbf{C}}=\textbf{T}_{i_{r}}^{\textbf{C}}\textbf{S}_{i,j}^{\textbf{C}}\textbf{T}_{j_c}^{\textbf{C}}$.

\subsection{Summary of overall algorithm at each non-leaf level}
The cluster bases products $\textbf{B}_{j}$ have been computed for all clusters $j$ before the MMP starts, since they are only related to the original cluster bases. 

At each non-leaf level, we do the following:
\begin{enumerate}
\item Collect four blocks in an $\textbf{NL}$ block in $\textbf{A}_{{\cal H}^2}$ to a block of $O(k_{l+1})$ size, using the newly generated children row cluster bases of $\textbf{C}$ (transfer matrices if children are not at the leaf level) and the original column cluster bases of $\textbf{B}$ (or transfer matrices). This is to generate the $(\textbf{NL}_{i,j}^{\textbf{A}})_{coll.}$, shown in (\ref{case2_coll}).
\item Collect four blocks in an $\textbf{NL}$ block in $\textbf{B}_{{\cal H}^2}$ to a block of $O(k_{l+1})$ size, using the original row cluster bases of $\textbf{A}$ (or transfer matrices) and the new children column cluster bases of $\textbf{C}$ (transfer matrices if children are not at the leaf level). This is to generate $(\textbf{NL}_{j,k}^{\textbf{B}})_{coll.}$, shown in (\ref{case3_coll}).
\item Merge four blocks in an $\textbf{R}$ block in $\textbf{C}_{{\cal H}^2}$. This corresponds to the $(\textbf{NL}_{i,j}^{\textbf{A}}\textbf{NL}_{j,k}^{\textbf{B}})_{coll.}^{(l)}$ in (\ref{nonleafcoup_simple}).
\item Calculate new row and column transfer matrices of product matrix $\textbf{C}_{{\cal H}^2}$ at this level.
\item Prepare cluster bases projections $\textbf{P}_{i}^{\textbf{A}}$, $\textbf{P}_{k}^{\textbf{B}}$, and perform an $\textbf{NL}$ block collect shown in (\ref{NLcollect}); 
\item Perform four cases of multiplications shown in (\ref{nonleafcoup_simple}).
\end{enumerate}

After we finish one-way bottom-up tree traversal to calculate block matrix products at all the levels, i.e. from leaf level all the way up to minimal admissible level, we need to perform a post-processing for the coupling matrices associated with the \textbf{NL} blocks in $\textbf{C}_{{\cal H}^2}$. They exist because of the multiplications cases described in Section \ref{sec:nl-nl}. This could be efficiently done by performing one-way top-down split process, the same as the matrix backward transformation shown in \cite{BormH2MMP}. This post processing stage is to split the coupling matrices in \textbf{NL} to lower level admissible or inadmissible blocks.

\section{Accuracy and Complexity Analysis}
In this section, we analyze the accuracy and computational complexity of the proposed algorithm to compute ${\cal H}^2$-matrix-matrix products.
\subsection{Accuracy}
Different from existing formatted ${\cal H}^2$-matrix-matrix products \cite{BormH2MMP}, in the proposed new algorithm, the accuracy of the product is directly controlled by $\epsilon_{trunc}$. No formatted multiplications are performed, and the cluster bases are changed to represent the updates to the original matrix accurately. This makes each operation performed in the proposed MMP controlled by accuracy or exact. When generating an ${\cal H}^2$-matrix to represent the original dense matrix, the accuracy is controlled by $\epsilon_{{\cal H}^2}$, which is the same as in \cite{SaadAP2015}. 

\subsection{Time and Memory Complexity}
The proposed MMP involves $O(L)$ levels of computation. At each level, there are $2^l$ clusters. For each cluster, the cost of changing the cluster bases at the leaf level due to four cases of multiplications is to perform $O(C_{sp})^2$ multiplications, and each of which has a constant cost, as can be seen from (\ref{case1-ch}), (\ref{case2-ch}), and (\ref{case3-ch}). The cost of changing the cluster bases at the non-leaf level due to the four cases of multiplications is also to perform $O(C_{sp})^2$ multiplications for each cluster, and each of which has a cost of $O(k_l)^3$, as can be seen from (\ref{case1-ch-nl}), (\ref{case2-ch-nl}), (\ref{case3-ch-nl}). Notice that the $\textbf{NL}$ blocks in \textbf{A} and \textbf{B} are collected level by level, at each level, there are $2^lO(C_{sp})$ $\textbf{NL}$ block, and each collect operation also costs $O(k_l)^3$ only. Other auxiliary matrices are generated using a similar computational cost.

As for the computation of the four cases of multiplications at each level, each case involves $O(C_{sp})^2$ multiplications for each cluster, and each of which costs $O(k_l)^3$ at the non-leaf level and ${O(leafsize)}^3$ at the leaf level as can be seen from (\ref{leafcoup_simple}), and (\ref{nonleafcoup_simple}).

Hence, the time complexity of the proposed MMP can be found as
\begin{equation}\label{timeComp}
\textbf{Time Complexity}=\sum_{l=0}^{L} C_{sp}^2 2^{l} O(k_l)^3= {C_{sp}^2}\sum_{l=0}^{L} 2^{l}O(k_l)^3.
\end{equation}
And the storage for each block is $O(k_{l}^2)$, with each cluster having $C_{sp}$ blocks. So the memory complexity is
\begin{equation}\label{memComp}
\textbf{Memory Complexity}=\sum_{l=0}^{L} C_{sp} 2^{l} O(k_l)^2= C_{sp}\sum_{l=0}^{L} 2^{l} O(k_l)^2.
\end{equation}
Recall $k_l$ is the rank at tree level $l$. Hence, (\ref{timeComp}) and (\ref{memComp}) show that the overall complexity is a function of rank $k_l$. Taking into account the rank's growth with electrical size as shown in \cite{JiaoRank}, we can get the time and memory complexity of proposed MMP for different rank scaling. For constant-rank ${\cal H}^2$-matrices, since $k_l$ is a constant irrespective of matrix size, the complexity of the proposed direct solution is strictly $O(N)$ in both CPU time and memory consumption, as shown below.

\textbf{For constant $k_l$:}
\begin{equation}\label{timeComp_const}
\textbf{Time Complexity}=C_{sp}^2 k_l^3 \sum_{l=0}^{L} 2^{l} =O(N),
\end{equation}
\begin{equation}\label{memComp_const}
\textbf{Memory Complexity}=C_{sp}k_l^2\sum_{l=0}^{L} 2^{l}=O(N).
\end{equation}
For electrodynamic analysis, to ensure a prescribed accuracy, the rank becomes a function of electrical size, and thereby tree level. Different ${\cal H}^2$-matrix representations can result in different complexities, because their rank's behavior is different. Using a minimal-rank ${\cal H}^2$-representation, as shown by \cite{JiaoRank}, the rank grows linearly with electrical size for general 3-D problems. In a VIE, $k_l$ is proportional to the cubic root of matrix size at level $l$, because this is the electrical size at level $l$. Hence for a VIE, (\ref{timeComp}) and (\ref{memComp}) become

\textbf{For $k_l$ linearly growing with electrical size:}
\begin{equation}\label{timeComp_linear}
\textbf{Time Complexity}={C_{sp}}^{2}\sum_{l=0}^{L} 2^{l} \left[ \left(\frac{N}{2^l}\right)^{\frac{1}{3}}\right]^3=O(NlogN),
\end{equation}
\begin{equation}\label{memComp_linear}
\textbf{Memory Complexity}=C_{sp}\sum_{l=0}^{L} 2^{l} \left[ \left(\frac{N}{2^l}\right)^{\frac{1}{3}} \right]^2=O(N).
\end{equation}
So the time complexity of the proposed MMP algorithm for 3D electrodynamic analysis is $O(NlogN)$, and the memory complexity is $O(N)$. 
\begin{figure*}[tbp] 
\vspace{-0.2in}
\centering
\subfloat[]{\includegraphics[width=0.49\textwidth]{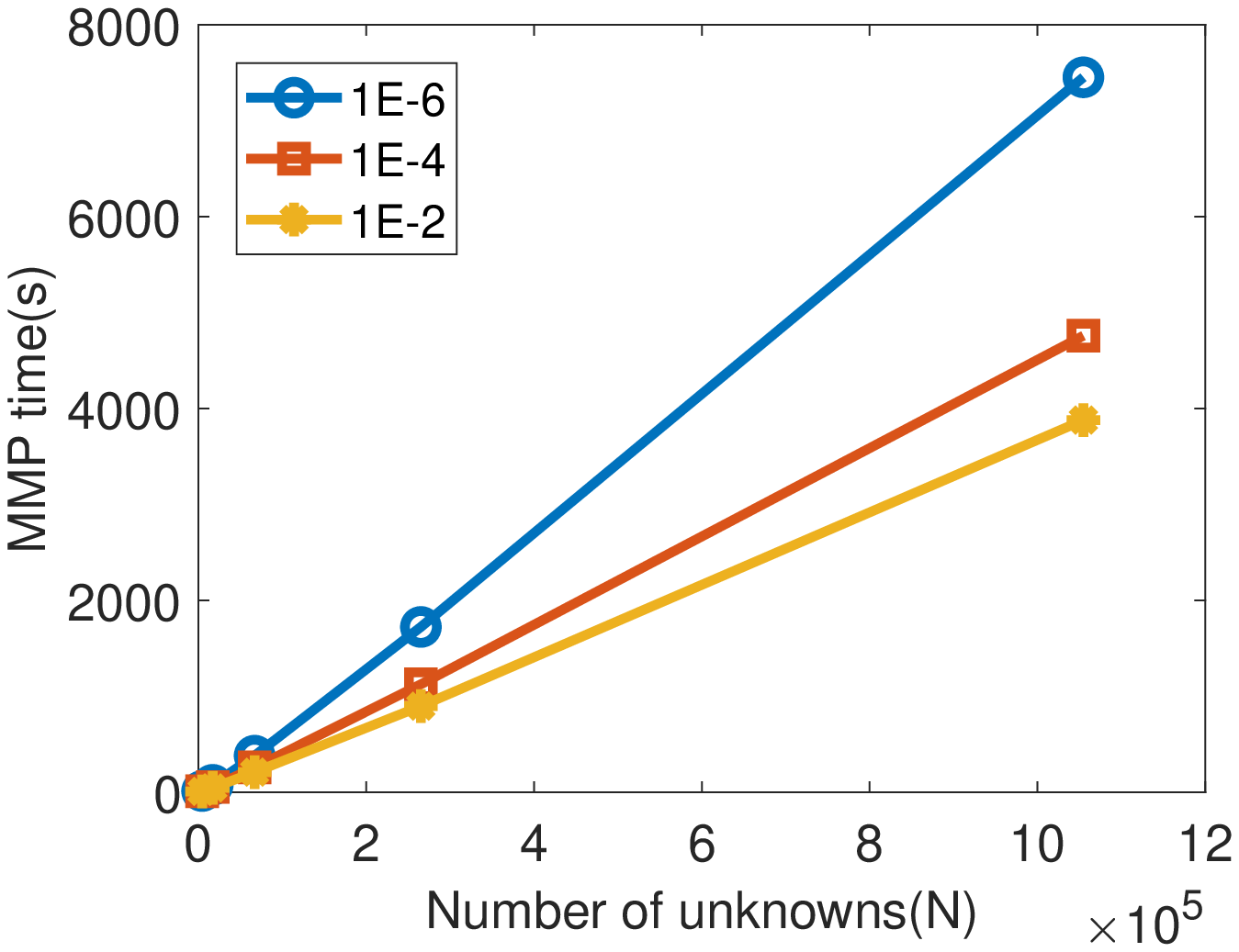}\label{fig:captime}}
\subfloat[]{\includegraphics[width=0.49\textwidth]{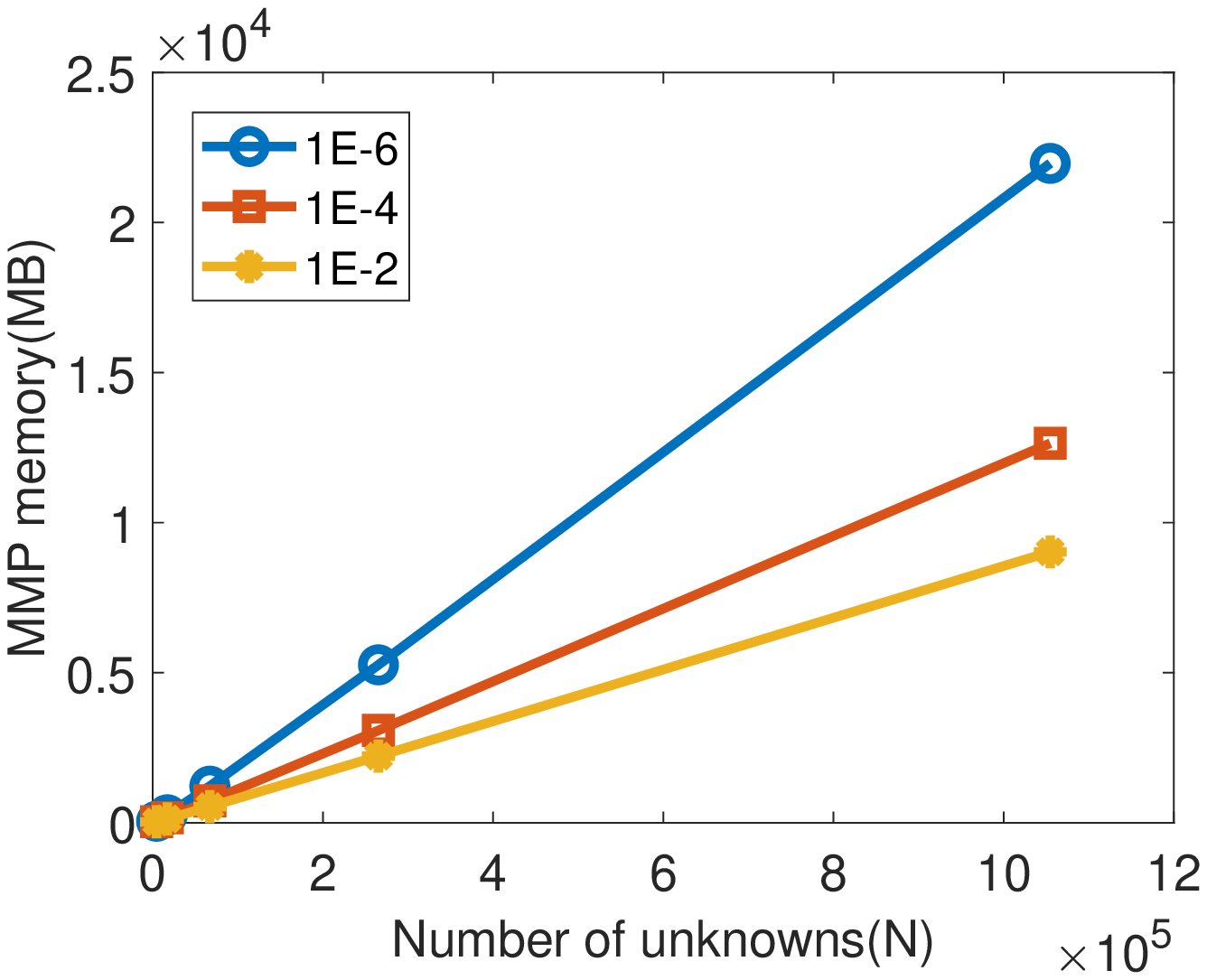}\label{fig:capmem}}
\caption{MMP performance for $\textbf{\textbf{A}}_{{\cal H}^2} \times \textbf{B}_{{\cal H}^2}$ of large-scale capacitance extraction matrices. (a) Time scaling v.s. $N$. (b) Memory scaling v.s. $N$.}
\label{captimeandmem}
\end{figure*}
\section{Numerical Results}
In order to demonstrate the accuracy and low computational complexity of the proposed fast ${\cal H}^2$-matrix-matrix multiplication for general ${\cal H}^2$-matrices, we use ${\cal H}^2$-matrices resulting from large-scale capacitance extraction and volume integral equation (VIE) solvers for electromagnetic analysis as examples. The capacitance extraction matrix is shown in \cite{Weninverse}. The VIE formulation is based on \cite{WiltonVIE} with SWG vector bases for expanding electric flux density in each tetrahedral element. A variety of large-scale examples involving over one million unknowns are simulated on a single CPU core to examine the accuracy and complexity of the proposed MMP algorithm. The ${\cal H}^2$-matrix for each example is constructed based on the method described in \cite{SaadICEEA2014, SaadAP2015}. The capacitance matrix is used to demonstrate the proposed MMP algorithm performance for constant-rank ${\cal H}^2$-matrices. We also simulate large scale 2- and 3-D scattering examples to examine the time and memory complexity of the proposed MMP for variable rank cases. The computer used has an Intel(R) Xeon(R) CPU E5-2690 v2 running at 3 GHz, and only a \emph{single core} is employed to carry out the computation. 

\subsection{Two-layer Cross Bus}
The first example is the capacitance extraction of a 2-layer cross bus structure. In each layer, there are $m$ conductors, and each conductor has a dimension of $1 \times 1 \times (2m+1)$ m$^3$. We simulate a suite of such structures with 16, 32, 64, 128, and 256 buses respectively. The parameters used in the ${\cal H}^2$-matrix construction are $leafsize = 30$, admissibility condition \cite{BormH2MMP} $\eta = 1.0$, and $\epsilon_{{\cal H}^2}=10^{-4}$. For the proposed ${\cal H}^2$ MMP, the $\epsilon_{trunc}$ is chosen to be $10^{-2}$, $10^{-4}$ and $10^{-6}$ respectively to examine the error controllability. As shown in Fig. \ref{captimeandmem}, the proposed MMP exhibit clear linear complexities in time and memory regardless of the choice of $\epsilon_{trunc}$. Certainly, the smaller the $\epsilon_{trunc}$, the larger the computational cost. 

The accuracy of the proposed MMP is assessed by using the following criterion:
\begin{equation}\label{error}
\epsilon_{rel}=\frac{||\textbf{C}_{{\cal H}^2} x-\textbf{A}_{{\cal H}^2} (\textbf{B}_{{\cal H}^2} x)||_{F}}{||\textbf{A}_{{\cal H}^2} (\textbf{B}_{{\cal H}^2} x)||_{F}},
\end{equation}
where $\textbf{A}_{{\cal H}^2} \times (\textbf{B}_{{\cal H}^2} \times x)$ is used as the reference solution, since given an ${\cal H}^2$ matrix, a matrix-vector product can be carried out without any approximation as shown in \cite{BormH2MMP}. In generating the reference solution, we first compute $y=\textbf{B}_{{\cal H}^2} \times x$, and then compute $\textbf{A}_{{\cal H}^2} \times y$, both of which are done in exact arithmetic. The proposed solution is generated by first computing an MMP of $\textbf{A}_{{\cal H}^2}\textbf{B}_{{\cal H}^2}$ to obtain $\textbf{C}_{{\cal H}^2}$, and then compute $\textbf{C}_{{\cal H}^2} x$. From Table \ref{tab:relErr_bus}, we can see the accuracy of the proposed MMP is good, and it is also controllable.
\begin{table}[h]
\caption{${\cal H}^2$ MMP error at different $\epsilon_{trunc}$ for large-scale capacitance extraction matrices as a function of $N$. \label{tab:relErr_bus}}
\centering
\small
\begin{tabular}{|c|c|c|c|c|c|}
\hline
$N$ &  4,480 & 17,152 & 67,072 & 265,216 & 1,054,720\\ \hline
$\epsilon_{tr}:$ 1E-2 & 4.35E-2 & 5.72E-2 & 5.73E-2 & 5.80E-2 & 5.97E-2 \\ \hline
$\epsilon_{tr}:$ 1E-4 & 3.71E-3 & 3.72E-3 & 3.86E-3 & 3.80E-3 & 3.67E-3 \\ \hline
$\epsilon_{tr}:$ 1E-6 & 2.82E-4 & 3.28E-4 & 3.86E-4 & 4.50E-4 & 5.66E-4 \\ \hline
\end{tabular}
\end{table}

\begin{figure*}[t]
\vspace{-0.2in}
\centering
\subfloat[]{\includegraphics[width=0.49\textwidth]{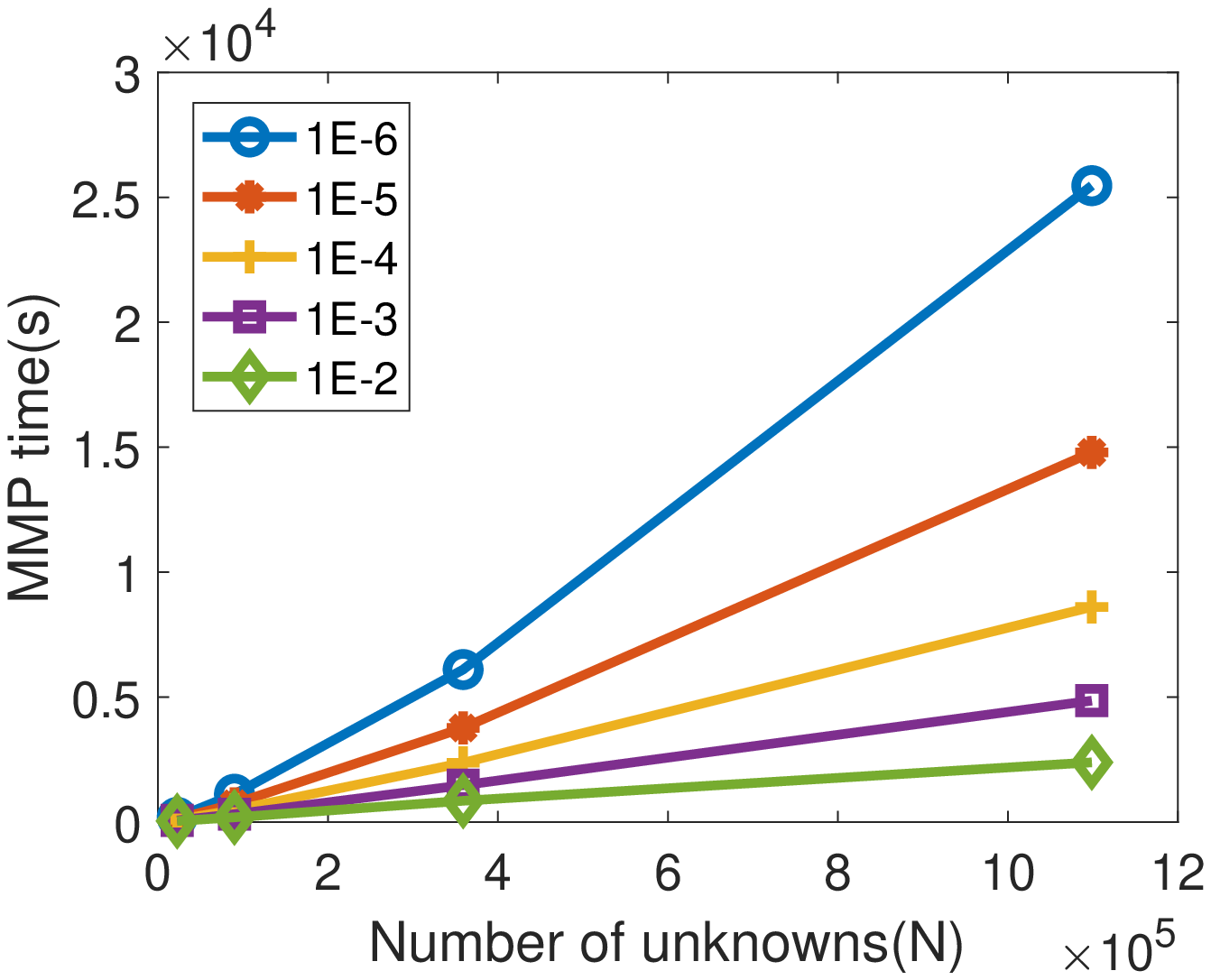}}
\subfloat[]{\includegraphics[width=0.49\textwidth]{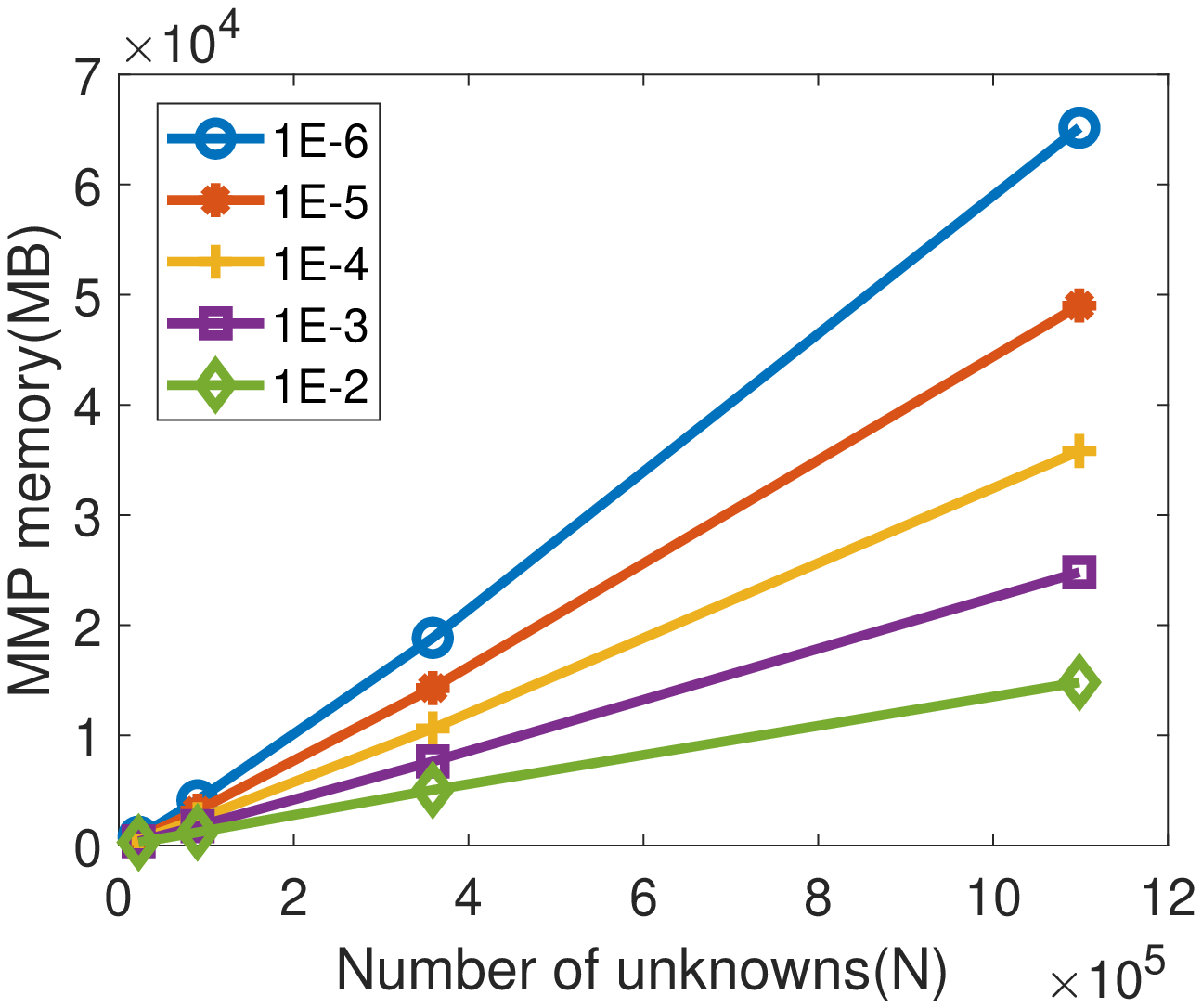}}
\caption{MMP performance for $\textbf{A}_{{\cal H}^2} \times \textbf{A}_{{\cal H}^2}$ of 2-D slab from $4\lambda$ to $28\lambda$. (a) Time scaling v.s. $N$. (b) Memory scaling v.s. $N$.}
\label{slabtimeandmem}
\end{figure*}
\begin{figure*}[t]
\centering
\subfloat[]{\includegraphics[width=0.49\textwidth]{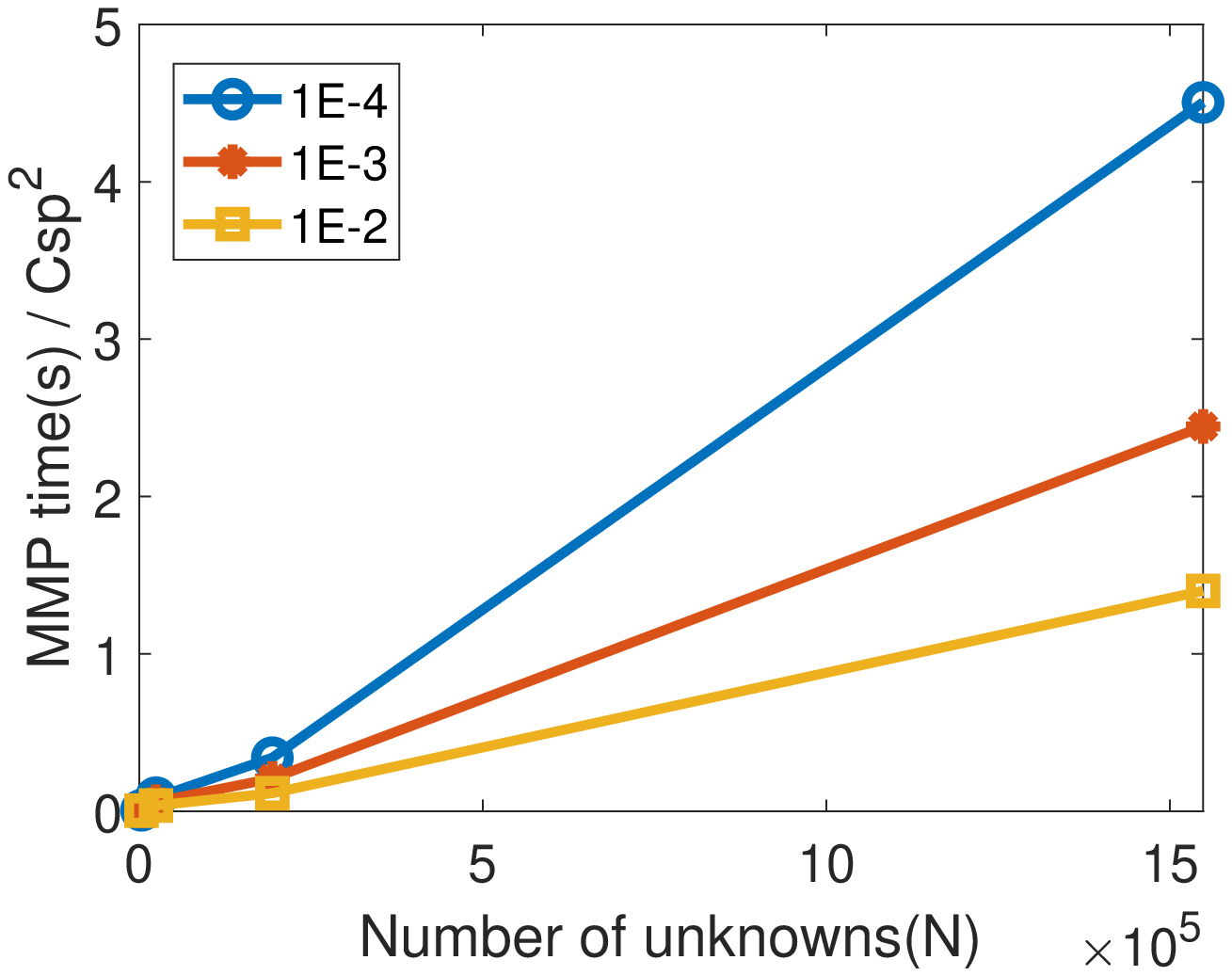}}
\subfloat[]{\includegraphics[width=0.49\textwidth]{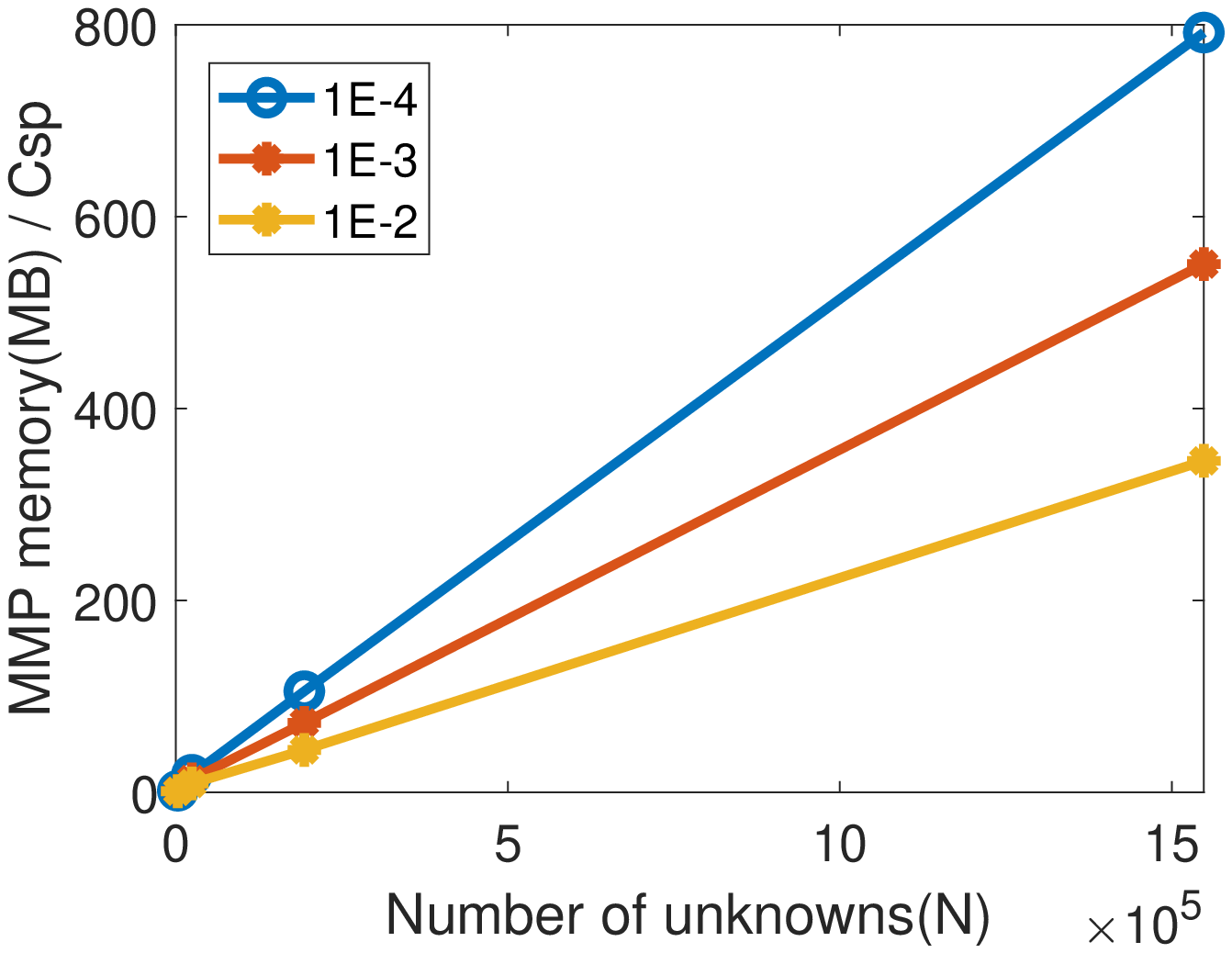}}
\caption{MMP performance for $\textbf{A}_{{\cal H}^2} \times \textbf{A}_{{\cal H}^2}$ of 3-D cube array. (a) Time scaling v.s. $N$. (b) Memory scaling v.s. $N$.}
\label{cubetimemem}
\end{figure*}

\subsection{Large-scale Dielectric Slab Scattering}
We then simulate a dielectric slab with $\epsilon_{r}=2.54$ at 300 MHz. The thickness of the slab is fixed to be  $0.1\lambda_0$. The width and length are simultaneously increased from $4\lambda_0$, $8\lambda_0$, $16\lambda_0$, to $28\lambda_0$. With a mesh size of $0.1\lambda_0$, the resultant $N$ ranges from 22,560 to 1,098,720 for this suite of slab structures. The parameters used in the ${\cal H}^2$-matrix construction are $leafsize = 40$, admissibility condition  \cite{BormH2MMP} $\eta = 2.0$, and $\epsilon_{{\cal H}^2}=10^{-3}$. For the proposed ${\cal H}^2$ MMP, the $\epsilon_{trunc}$ is chosen to be $10^{-2}$, $10^{-3}$, $10^{-4}$, $10^{-5}$and $10^{-6}$ respectively, to examine the computational complexity and error controllability of the proposed MMP. Based on \cite{JiaoRank}, the rank's growth rate with electrical size for 2-D slab is lower than linear, and being a square root of the log-linear of the electric size. Substituting such a rank's growth into the complexity analysis in (\ref{timeComp}) and (\ref{memComp}), we obtain linear complexity in both memory and time. 

In Fig. \ref{slabtimeandmem} (a), we plot the MMP time with respect to $N$, for all different choices of $\epsilon_{trunc}$. It is clear that the smaller $\epsilon_{trunc}$ value, the larger the MMP time. However, the complexity remains the same as linear regardless of the choice of $\epsilon_{trunc}$. The memory cost is plotted in Fig. \ref{slabtimeandmem} (b). Obviously, it scales linearly with the number of unknowns. The error of the proposed MMP is measured in the same way as shown in (\ref{error}). In Table \ref{slab_acc}, we list the error as a function of $\epsilon_{trunc}$. Excellent accuracy can be observed in the entire unknown range. Furthermore, the accuracy can be controlled by $\epsilon_{trunc}$, and overall smaller $\epsilon_{trunc}$ results in better accuracy.

\begin{table}[h]
\caption{${\cal H}^2$ MMP error for 2-D slab.\label{tab:relErr_slab} for different $\epsilon_{trunc}$ as a function of $N$.}
\label{slab_acc}
\centering
\begin{tabular}{|c|c|c|c|c|}
\hline
$N$ &  22560 & 89920 & 359040 & 1098720\\ \hline
1E-2 & 8.54E-3 & 1.06E-2 & 1.49E-2 & 1.07E-2  \\ \hline
1E-3 & 2.52E-3 & 3.17E-3 & 4.23E-3 & 3.79E-3 \\ \hline
1E-4 & 7.86E-4 & 9.76E-4 & 1.38E-3 & 1.23E-3  \\ \hline
1E-5 & 2.91E-4 & 3.37E-4 & 4.22E-4 & 4.11E-4  \\ \hline
1E-6 & 8.04E-5 & 9.85E-5 & 1.27E-4 & 1.36E-4  \\ \hline
\end{tabular}
\end{table}
\subsection{Scattering from Large-scale Array of Dielectric Cubes}
Next, we simulate a large-scale array of dielectric cubes at 300 MHz. The relative permittivity of the cube is $\epsilon_{r}=4.0$. Each cube is of size $0.3\lambda_0 \times 0.3\lambda_0 \times 0.3\lambda_0$. The distance between adjacent cubes is kept to be $0.3\lambda_0$. The number of the cubes is increased along the $x$-, $y$-, and $z$- directions simultaneously from 2 to 16, thus producing a 3-D cube array from $2 \times 2 \times 2$ to $16 \times 16 \times 16$ elements. The number of unknowns $N$ is respectively 3,024, 24,192, 193,536, and 1,548,288 for these arrays. During the construction of ${\cal H}^2$-matrix, we set $leafsize = 20$, $\eta =1$ and $\epsilon_{{\cal H}^2}=10^{-2}$. For the proposed ${\cal H}^2$ MMP, the $\epsilon_{trunc}$ is chosen as $10^{-2}$, $10^{-3}$ and $10^{-4}$.
\begin{table}[htpb]
\caption{$C_{sp}$ as a function of $N$ for the dieletric cube array}
\label{csp}
\centering
\begin{tabular}{|c|c|c|c|c|}
\hline
$N$ &  3024 & 24192 & 193536 & 1548288\\ \hline
$C_{sp}$ & 16 & 42 & 95 & 126  \\ \hline
\end{tabular}
\end{table}

For a cubic growth of unknowns in 3-D problems, we observe that constant $C_{sp}$ is quite different for different unknowns, as can be seen from Table \ref{csp}. It is thus important to analyze the performances of the proposed MMP as $Memory/C_{sp}$ and $Multiplication\ time/C_{sp}^{2}$ respectively to examine the true scaling rate. In Fig. \ref{cubetimemem} (a) and Fig. \ref{cubetimemem} (b), we plot the ${\cal H}^2$-matrix-matrix multiplication time divided by $C_{sp}^{2}$, and the storage cost normalized with $C_{sp}$ with respect to N. As can be seen, their scaling rate with $N$ agrees very well with our theoretical complexity analysis. For the over one-million unknown case which is a $16\times16\times16$ cube array
having thousands of cube elements, the error is still controlled
to be as small as $0.809\%$ using $\epsilon_{trunc}=10^{-4}$. The error of the proposed MMP is listed in Table \ref{tab:relErr_cube} for this example, which again reveals excellent accuracy and error controllability of the proposed MMP. 
\begin{table}[htpb]
\caption{${\cal H}^2$ MMP error at different $\epsilon_{trunc}$ for 3-D cube array.\label{tab:relErr_cube}}
\centering
\begin{tabular}{|c|c|c|c|c|}
\hline
$N$ &  3024 & 24192 & 193536 & 1548288\\ \hline
Existing \cite{BormH2MMP} & 9.02E-2 & 1.01E-1 & 1.77E-1 & 2.74E-1 \\ \hline
1E-2 & 1.91E-2 & 2.38E-2 & 3.82E-2 & 6.58E-2 \\ \hline
1E-3 & 5.51E-3 & 7.23E-3 & 1.06E-2 & 2.16E-2  \\ \hline
1E-4 & 1.48E-3 & 2.46E-3 & 3.69E-3 & 8.09E-3 \\ \hline
\end{tabular}
\end{table}
We also compare the accuracy of the proposed MMP with existing MMP \cite{BormH2MMP} using this 3-D example. As shown in Table \ref{tab:relErr_cube}, the proposed MMP has much better accuracy, and also it is controllable.
% \begin{figure}[h] 
% \centering
% \subfloat[]{\includegraphics[width=0.49\textwidth]{Cube_AA_Compare.eps}}
% \caption{MMP time comparison with existing MMP for 3D cube array VIE matrix.}
% \label{3Dcompare}
% \end{figure}
% \begin{table}[htpb]
% \caption{${\cal H}^2$ MMP error comparison with existing MMP for 3-D cube array.\label{tab:cubecomp}}
% \centering
% \small
% \begin{tabular}{|c|c|c|c|c|}
% \hline
% $N$ &  3024 & 24192 & 193536 & 1548288\\ \hline
% Existing \cite{BormH2MMP} & 9.02E-2 & 1.01E-1 & 1.77E-1 & 2.74E-1  \\ \hline
% Proposed & 8.03E-2 & 7.39E-2 & 1.03E-1 & 1.25E-1  \\ \hline
% \end{tabular}
% \end{table}

\section{Conclusions}
In this paper, we develop a fast accuracy-controlled algorithm to compute ${\cal H}^2$-matrix-matrix products for general ${\cal H}^2$-matrices. This proposed algorithm not only has \emph{explicitly} controlled accuracy, but also generates a rank-minimized representation of the product matrix based on prescribed accuracy. The row and column cluster bases are instantaneously changed so that the new matrix content generated during the MMP can be accurately represented. This ensures that each multiplication performed in the proposed MMP is well controlled by accuracy. Meanwhile, we retain the complexity to be linear for constant-rank ${\cal H}^2$-matrices. The proposed algorithm has been applied to calculate ${\cal H}^2$-matrix-matrix products for large-scale capacitance extraction matrices whose kernel is static and real-valued and electrically large VIEs whose kernel is oscillatory and complex-valued. For constant-rank ${\cal H}^2$-matrices, the proposed MMP has an $O(N)$ complexity in both time and memory. For rank growing with the electrical size linearly, the proposed MMP has an $O(NlogN)$ complexity time and $O(N)$ complexity in memory. ${\cal H}^2$-matrix products with millions of unknowns are simulated on a single core CPU in fast CPU run time. Comparisons with existing ${\cal H}^2$-matrix-matrix product algorithm have demonstrated clear advantages of the proposed new MMP algorithm.

% or


%\appendices

%\section{}

%\section*{Acknowledgment}


\ifCLASSOPTIONcaptionsoff
  \newpage
\fi

\end{document}